\documentstyle[11pt,leqno]{article}
\input amssym.def
\input amssym.tex
\def\CC{{\Bbb C}}

\def\NN{{\Bbb N}}

\def\RR{{\Bbb R}}

\def\l{\langle}
\def\ra{\rangle}

\def\A{{\cal A}}
\def\B{{\cal B}}
\def\D{{\cal D}}
\def\F{{\cal F}}
\def\E{{\cal E}}
\def\H{{\cal H}}
\def\K{{\cal K}}

\def\M{{\cal M}}
\def\cS{{\cal S}}

\def\U{{\cal U}}

\def\Lp{{L^{\u p}(G/P,\ell)}}
\def\Lq{{L^{\u q}(G/P,\ell)}}

\def\a{{\frak a}}
\def\b{{\frak b}}
\def\g{{\frak g}}
\def\h{{\frak h}}
\def\m{{\frak m}}
\def\n{{\frak n}}
\def\p{{\frak p}}
\def\r{{\frak r}}
\def\s{{\frak s}}
\def\uu{{\frak u}}
\def\vv{{\frak v}}
\def\ww{{\frak w}}
\def\k{{\frak k}}

\def\al{\alpha}
\def\de{\delta}
\def\la{\lambda}

\def\ve{\varepsilon}

\def\vp{\varphi}
\def\Om{\Omega}
\def\om{\omega}
\def\si{\sigma}
\def\ga{\gamma}
\def\ze{\zeta}
\def\De{\Delta}

\def\Si{\Sigma}

\def\o{\overline}
\def\u{\underline}

\def\ad{{\rm ad}}

\def\Ad{{\rm Ad}}

\def\ind{{\rm ind}}
\def\ker{{\rm ker}}

\def\Im {{\rm Im}\,}
\def\Re {{\rm Re}\,}
\def\span{{\rm span}}

\def\tr{{\rm tr}}
\def\Prim{{\rm Prim}}

\def\Id{{\Bbb I}}
\def\im{{\rm im}\,}
\def\inv{^{-1}}

\def\be{\begin{enumerate}}
\def\ee{\end{enumerate}}
\def\noi{\noindent}

\newtheorem{theorem}[subsection]{Theorem}
\newtheorem{prop}[subsection]{Proposition}
\newtheorem{cor}[subsection]{Corollary}
\newtheorem{lemma}[subsection]{Lemma}

\begin{document}

\title{Sub-Laplacians of holomorphic $L^p$-type
on exponential solvable groups}

\author{ W. Hebisch, J. Ludwig  and D. M\"uller}

\date{}

\maketitle

\begin{abstract}Let $L$ denote a right-invariant sub-Laplacian on an
exponential,
hence solvable Lie group~$G$, endowed with a left-invariant Haar measure.
Depending
 on the structure of $G$, and possibly also that of $L$, $L$ may admit
differentiable $L^p$-functional calculi, or may be of holomorphic
$L^p$-type for
a given $p\ne 2$. By ``holomorphic $L^p$-type'' we mean that every
$L^p$-spectral
multiplier for $L$ is necessarily holomorphic in a complex neighborhood of
some non-isolated point of the $L^2$-spectrum of $L$. This can in fact only
arise
if the group algebra $L^1(G)$ is non-symmetric.

Assume that $p\ne 2$. For a point $\ell$ in the dual  $\frak g ^*$  of the Lie
algebra $\frak g$ of $G$, we denote by $\Omega(\ell)=Ad^*(G)\ell$  the
corresponding
 coadjoint orbit.  We prove that every sub-Laplacian on $G$ is of holomorphic
$L^p$-type, provided there exists a point $\ell\in \frak g ^*$   satisfying
``Boidol's condition'' (which, by \cite{Poguntke} is equivalent to the
non-symmetry
of $L^1(G)$), such that  the restriction of $\Omega(\ell)$ to the
nilradical of $\frak g$
is closed.
  This work improves on the results in
\cite{ChristMueller}, \cite{Ludwig-Mueller} in twofold ways: On the
one hand, we no longer impose any restriction on the structure of the
exponential
group $G$, and on the other hand, for the case $p>1$,  our conditions need
to hold
for a single coadjoint orbit only, and not for an open set of orbits, in
contrast
to \cite{Ludwig-Mueller}.

It seems likely that the condition that the restriction of $\Omega(\ell)$ to
the nilradical of $\frak g$ is closed could be replaced by the weaker
condition that
the orbit $\Omega(\ell)$ itself is closed. This would then prove one
implication
of a conjecture  made in \cite{Ludwig-Mueller}, according to which there exists a sub-Laplacian of holomorphic $L^1$ (or, more generally, $L^p$)-type on
$G$ if and
only if there exists a point $\ell\in \frak g ^*$ whose orbit  is closed
and which
satisfies Boidol's condition.
\hspace{-0.2cm}\footnote{  keywords: exponential solvable Lie group,
sub-Laplacian, functional
calculus, $L^p$-spectral multiplier,  symmetry, unitary representation,
holomorphic $L^p$-type, heat kernel, transference, perturbation }
\end{abstract}

\section*{Introduction}

A comprehensive discussion of the problem  studied in this article, background
information and references to further literature  can be found in
\cite{Ludwig-Mueller}.
 We shall therefore content ourselves in this introduction by recalling some
 notation and results from \cite{Ludwig-Mueller}.

If $T$ is a self-adjoint linear operator on a Hilbertian $L^2$-space
 $L^2(X,d\mu)$, with spectral resolution $T=\int\limits_{\RR}
 \lambda dE_\la$, and if $m$ is a bounded Borel function on $\RR$, then we
call
$m$ an $L^p$-{\em multiplier for } $T$ ($1\le p<\infty$), if
 $m(T):=\int_\RR m(\la)dE_\la$ extends from $L^p\cap L^2(X,d\mu)$
to a bounded operator on $L^p(X,d\mu)$. We shall denote by $\M_p(T)$
  the space of all $L^p$-multipliers for $T$, and by $\si_p(T)$ the $L^p$-
spectrum of $T.$  We say
that $T$ is of {\em holomorphic $L^p$-type}, if there exist some
non-isolated point $\la_0$ in the $L^2$-spectrum $\si_2(T)$ and
 an open complex neighborhood $\U$ of $\la_0$ in $\CC$, such that
 every $m\in \M_p(T)\cap C_\infty (\RR)$ extends holomorphically to
 $\U$. Here, $C_\infty (\RR)$ denotes the space of all continuous
functions on $\RR$ vanishing at infinity.

Assume in addition that there exists a linear subspace $\D$ of $L^2(X)$
which is $T$-invariant and  dense in $L^p(X)$ for every $p\in [1,\infty[$,
and that
$T$ coincides with the closure of its restriction to $\D$. Then, if $T$ is
of  holomorphic $L^p$-type, the set $\U$
 belongs to the $L^p$-spectrum of $T$, i.e.
\begin{equation}\overline\U\subset\si_p(T).
\end{equation}
In particular,
\begin{equation}\si_2(T)\subsetneq \si_p(T).
\end{equation}

Throughout this article, $G$ will denote an {\em exponential} Lie group, i.e.
the exponential mapping $\exp:\g\to G$ is a diffeomorphism from the Lie
algebra
$\g$ of $G$ onto $G$. Such a group is solvable \cite{Bernat}.
The  inverse mapping to $\exp$ will be denoted by $\log$.

We fix a {\em left-invariant Haar measure} $dg$ on $G$.
If $\pi:G\to \cal U(\cal H)$ is a unitary representation of $G$ on the Hilbert
space $\H=\H_\pi$, then we denote the integrated representation of
$L^1(G)=L^1(G,dg)$ again by $\pi$, i.e. $\pi(f)\xi:=\int _G
f(g)\pi(g)\xi\,dg$ for
 every $f\in L^1(G), \ \xi\in\H.$  For $X\in \g$, we denote by $d\pi(X)$
the infinitesimal
 generator of the one-parameter group of unitary operators
$t\mapsto \pi(\exp tX)$. By $X^r$  we denote the
right- invariant vector field on $G$, given by
\begin{eqnarray*}
X^r f(g)&:=& \lim_{t\to 0} \frac 1 t [f((\exp tX)g)-f(g)].
\end{eqnarray*}

\noi For a given function $f$ on $G$, we write
\[
[\la(g)f](x):=f(g^{-1}x),  \quad  g,x\in G,
\]
for the {\em left-regular}  action of $G$. Then $\la$, acting on $L^2(G)$,
is a
unitary representation. In particular, we have
\begin{equation}
X^r=-d\la(X)
\end{equation}
and

\begin{equation}
\pi(X^r\vp)=-d\pi(X)\pi(\vp)
\end{equation}
for every $X\in \g,\  \vp\in \D(G):=C_0^\infty (G)$ and every unitary
representation
$\pi$ of $G$.
\medskip

 In the sequel, we shall usually identify $X\in \g$ with
the right-invariant vector field $-X^r=d\la(X)$, since $d\la$ (as $d\pi$
for any
unitary representation $\pi$) is a morphism of Lie algebras. One should notice
that $d\la(X)$ agrees with $-X$ at the identity $e$ of $G$, not with $X$. Then
(0.4) reads simply

\begin{equation}
\pi(X\vp)=d\pi(X)\pi(\vp).
\end{equation}

$d\pi$ extends from $\g$ to a representation $\pi_\infty$ of the universal
 enveloping algebra $\uu (\g)$ of $\g$ on the space $C^\infty(\pi)$ of all
$C^\infty$-vectors for $\pi$. Extending the convention above, we shall
often identify $A\in \uu(\g)$ with the right-invariant differential operator
 $\la_\infty(A)$ on $G$. Notice that $\la_\infty (\uu(\g))$ consists of all
 right-invariant complex coefficient differential operators on $G$.

Choose right invariant vector fields $X_1,\dots,X_k$ of $\g$   generating
$\g$ as a Lie algebra, and form the so-called
 {\em sub-Laplacian}
$$L=-\sum_{j=1}^k X_j^2.$$
By \cite{Nelson}, \cite{Hoermander} $L$ is hypoelliptic and essentially
self-adjoint as an operator on $L^2(G,dg)$ with domain $\D (G)$.
 We denote its closure again by L. Since $G$ is amenable, one has
\begin{equation}\si_2(L)=[0,\infty[.
\end{equation}

In this article, we shall give sufficient conditions for such an operator to
be of holomorphic $L^p$-type. As has been explained in \cite{Ludwig-Mueller},
a necessary condition for this to happen is the non-symmetry of the underlying
group. Recall that the {\em modular function} $\Delta_G$ on $G$  is defined by
the equation

$$
\int_G f(xg)dx=\Delta_G (g)^{-1} \int_Gf(x)dx, \quad g\in G.
$$

We put
\begin{eqnarray*}
\check f(g)&:=& f(g^{-1}),\nonumber\\
f^*(g)&:=& \Delta_G^{-1}(g) \o {f(g^{-1})}.
\end{eqnarray*}
Then $f\mapsto f^*$ is an isometric involution on $L^1(G)$, and
 for any unitary representation $\pi$ of $G$, we have
 $$
\pi(f)^*=\pi(f^*)\ .
$$

The group $G$ is said to be {\em symmetric}, if the associated group algebra
$L^1(G)$ is symmetric, i.e. if every element $f\in L^1(g)$ with $f^*=f$ has
 a real spectrum with respect to the involutive Banach algebra $L^1(G).$

The exponential  solvable non-symmetric Lie groups have been completely
classified by Poguntke \cite{Poguntke} (with previous contributions
by Leptin, Ludwig and Boidol)  in terms
of a purely Lie-algebraic condition (B). Let us describe this condition,
 which had been first introduced by Boidol in a different
context \cite{Boidol}.

Recall that the unitary dual of $G$ is in one to one correspondence with
the space of coadjoint orbits in $\g*$ via the Kirillov map, which
associates with a given point $\ell\in\g*$ an irreducible unitary
representation
$\pi_{\ell}$ (see Section 1).

If $\ell$ is an element of the dual space $\g^*$ of $\g$, denote by
 $$\g(\ell):=\ker\, \ad^*(\ell)=\{X\in \g:\ell([X,Y])=0 \ \forall Y\in \g\}$$
 the {\em stabilizer } of $\ell$ under the coadjoint
action $\ad^*$. Moreover, if $\m$ is any Lie algebra, denote by
$$\m=\m^1\supset \m^2\supset \dots$$
 the descending central series of $\m$, i.e. $\m^2=[\m,\m]$, and
$\m^{k+1}=[\m,\m^k]$.
 Put
$$\m^\infty =\bigcap_k\m^k.$$
 $\m^\infty$ is the smallest ideal $\k$ in $\m$ such that $\m/\k$ is
nilpotent. Put
 $$\m(\ell):=\g(\ell)+[\g,\g].$$
 Then we say that $\ell$ respectively the associated coadjoint orbit
$\Om(\ell):=\Ad^*(G)\ell$ satisfies {\em Boidol's condition}\,  (B), if

\noi
\parbox{1cm}{(B)} \hspace{4cm}\parbox{7cm}{ \begin{eqnarray*}
\ell\mid_{\m(\ell)^\infty} \ne 0.\end{eqnarray*}}
\medskip

According to \cite{Poguntke}, the group $G$ is non-symmetric if and only
if there exists a coadjoint orbit satisfying Boidol's condition.

If $\Om$ is a coadjoint orbit, and if $\n$ is the nilradical of $\g$, then
$$\Om|_{\n}:=\{\ell|_{\n}: \ell\in\Om\}\subset \n^*$$
will denote the restriction of $\Om$ to $\n.$

In this article, we shall prove the following extension and improvement
of the main theorems in \cite{Ludwig-Mueller}.

\bigskip

\noi {\bf Theorem 1} {\it  Let $G$ be an exponential solvable Lie group, and
assume that there exists a coadjoint orbit $\Om(\ell)$ satisfying Boidol's
condition, whose restriction
to the nilradical $\n$ is closed. Then every sub-Laplacian on $G$ is of
holomorphic
$L^p$-type, for $1\le p<\infty.$}
\bigskip

\noi {\bf Remarks.} (a) A sub-Laplacian $L$ on $G$ is of holomorphic $L^p$-type
 if and only if every continuous bounded multiplier $F\in\M_p(L)$ extends
holomorphically to an open neighborhood of a non-isolated point in $\si_2(L).$

For, with $F$, also the function $\tilde F(\la):=e^{-\la}F(\la)$ lies in
$\M_p(L)$,
since $\tilde F(L)=e^{-L}F(L),$ where the heat operator $e^{-L}$ is bounded
on every
$L^p(G).$ Furthermore, $\tilde F$ lies in $C_\infty(\RR).$

(b) If the restriction of a coadjoint orbit to the nilradical
is closed, then the orbit itself is closed (see Thm. 2.2).

(c) Under the hypotheses of the theorem, we obtain in particular that the
$L^2$-spectrum of $L$ is strictly contained in the $L^p$-spectrum of $L$ (see
(0.2)). This results has been proved independently by D. Poguntke \cite
{Poguntke-oral}.

(d)  What we really use in the proof is the following property of the orbit
$\Om:$
\medskip

\noi  {\em $\Om$ is closed, and for every real character $\nu$ of $\g$
which does
not vanish on $\g(\ell),$ there exists a sequence $\{\tau_n\}_n$ of real
numbers
such that $\lim_{n\to\infty}\Om+\tau_n\nu=\infty$ in the orbit space.}

\medskip
This property is a consequence of the closedness of $\Om|_\n.$
There are, however,  many examples where  the condition above is satisfied, so
that the conclusion of the theorem still holds, even though the
restriction of
$\Om$ to the nilradical is not closed (see e.g. Section 7). We do not know
whether the condition above automatically holds whenever the orbit $\Om$ is
 closed. 

\bigskip
The article is organized as follows: In Sections 1 and 2 we recall some
basic facts
from the unitary representation theory of exponential Lie groups (compare
\cite{Bernat},
\cite{LeptinLudwig}). Moreover, we prove a kind of Riemann-Lebesgue lemma
for one
parameter families of coadjoint orbits whose restrictions to the nilradical
are closed.
In the third section, we show how the irreducible unitary representations
of such
a group, which are in fact induced from characters of suitable polarizing
subgroups,
can be realized on Euclidean $L^2$-spaces. This will then allow for the
construction
 of analogous, isometric representations on certain mixed $L^p$-spaces.
Section 4
 provides some auxiliary results. In Section 4.1, we prove some results on
compact
operators acting on mixed $L^p$-spaces and their spectral properties. In
particular,
we prove an extension of a classical interpolation theorem by Krasnoselskii
for
compact operators acting on mixed $L^p$-spaces. Moreover, making use of
well-known
results on approximate units of Herz-Schur multipliers for amenable groups,
we prove
a result on the approximation of certain convolution operators by
convolutions with
continuous functions with compact support. This result will later allow us
to apply
a transference result by Coifman and Weiss to spectral multiplier operators
$F(L).$
 In Section 5 we show how, in the presence of Boidol's condition, one can
construct
certain analytic families $\{\pi_\ell^z\}_z$ of bounded representations
acting on
mixed $L^p$- spaces. Moreover, putting $T(z):=\pi_\ell^z(h_1),$ where $h_t$
denotes
the heat kernel associated to $L$ at time $t>0,$ we show that $\{T(z)\}_z$
is an analytic
family of compact operators on a wide range of mixed $L^p$-spaces, so that
we can
apply analytic perturbation theory.  Putting  together all results from the
preceding
sections, we complete the proof of Theorem 1 in Section 6. Finally, in
Section 7 we
present the  example announced in Remark (d).

\setcounter{equation}{0}
\section{Irreducible unitary representations}

Let again $G=\exp \g$ denote an exponential solvable Lie group and $\n$ a
nilpotent ideal of $\g$ containing $[\g,\g]$. Consider a composition sequence
\[
\g=\g_0\supset \g_1\supset\dots\supset \g_m=\{0\}
\]
for the adjoint action of $\g$, so that $\g_j/\g_{j+1}$ is an irreducible
 $\ad(\g)$-module. Since $\g$ is solvable, by Lie's theorem
we have $\dim \g_j/\g_{j+1}\le 2$. We may and shall assume that $\g_q=\n$ for
 some $q$. Choose a refinement
\[
\g=\a_0\supset \a_1\supset \dots\supset \a_r=\{0\}
\]
of the composition sequence, which means that $\dim(\a_j/ \a_{j+1})=1$, and
that either $\a_j=\g_i$ for some $i$, or, if $\a_j$ is not an ideal of $\g$,
then $\a_{j-1}=\g_i$ and $\a_{j+1}=\g_{i+1}$ for some $i$. We call such a
sequence $\{\a_j\}_j$ a {\em Jordan-H\"older sequence} for $\g$. Each $\a_j$ is
 a subalgebra of $\g$.

Let now $\ell$ be an element of $\g^*$. Denote by $\a_j(\ell)$ the subalgebra
 $\a_j(\ell):=\{X\in \a_j:\ell([X,\a_j])=\{0\}\}$, i.e. $\a_j(\ell)$ is the
stabilizer of $\ell|_{\a_j}$ in $\a_j$.

Put
\[
\p(\ell):=\sum_{j=0}^{r-1} \a_j(\ell).
\]
Then $\p(\ell)$ is a so-called {\em Vergne-polarization} for $\ell$. In
particular, it is a  {\em polarization}, i.e. a subalgebra $\p$ of $\g$ of
 maximal possible  dimension $\frac 1 2(\dim \g+\dim \g(\ell))$ such that
$\ell([\p,\p])=\{0\}$. Let $P( \ell):=\exp \p(\ell)\subset G$. We can define
the unitary character
\[
\chi_\ell(p):=e^{i\ell(\log p)} , \quad p\in P(\ell),
\]
of the closed subgroup $P(\ell)$, and denote by
\[
\pi_\ell=\pi_{\ell,P(\ell)}:=\ind _{P(\ell)}^G \chi_\ell
\]
the unitary representation of $G$ induced by the character $\chi_\ell$ of
$P(\ell)$. Let us briefly recall the notion of {\em induced representation }
 \cite{Bernat}:

If $P$ is any closed subgroup of $G$, with left-invariant Haar measure $dp$ and
 modular function $\Delta_P$, denote for $F\in C_0(G)$ by $\dot F$ the
function on
$G$ given by
\[
\dot F(x):=\int_PF(xp)\frac{\Delta _G}{\Delta_P} (p)\,dp,\quad x\in G,
\]
where $\Delta_G$ and $\Delta_P$ denote the modular functions of $G$ and $P,$
respectively. We shall also write
$\Delta_{G,P}$ instead of $\Delta_G/\Delta_P$. Then $\dot F$ lies in the space

\begin{eqnarray*}\E (G,P):=\{f\in C(G,\CC): f\ \mbox{ has compact support
modulo}\
 P,\\
 \ \mbox {and}\  f(xp)=(\Delta_{G,P}(p))^{-1}f(x)\quad \forall x\in G, p\in
P\}.
\end{eqnarray*}

In fact, one can show that $\E(G,P)=\{\dot F:F\in C_0(G)\}$. Moreover,
one checks that $\dot F=0$ implies $ \int_G F(x)dx=0$. From here it follows
that there exists
 a unique positive linear functional, denoted by $\oint_{G/P} \,d\dot x$,
on the
space $\E(G,P)$, which is left-invariant under $G$, such that

\begin{equation}
\int_G F(x)\,dx=\oint_{G/P} \dot F(x)\,d\dot x =\oint_{G/P} \int
f(xp)\Delta_{G/P}(p)\,dp\,d\dot x
\end{equation}
for every $F\in C_0(G)$.
\vskip.4cm

Now, given $\ell$ and the polarizing subgroup $P=P(\ell)$, put

\begin{eqnarray*}\E(G,P,\ell):=\{f\in C(G,\CC): f \ \mbox{ has compact
support modulo}
 P, \\
\ \mbox{and} \ f(xp)=\o{\chi_\ell(p)} (\Delta_{G,P}(p))^{-1/2} f(x)
\ \forall x\in G,& p\in P\},\end{eqnarray*}

\noindent  endowed with the norm
\[
||f||_2:= \left( \oint_{G/P} |f(x)|^2\,d\dot x\right)^{1/2}.
\]
Observe that $|f|^2\in \E(G,P)$. Let $\H_\ell=\H_{\ell,P(\ell)}$ denote
the completion of $\E(G,P,\ell)$ with respect to this norm. Then $\H_\ell$
becomes a Hilbert space, on which $G$ acts by left-translations  isometrically,
 and $\pi_\ell$ is defined  on $\H_\ell$ by
\[
[\pi_\ell(g)f](x):=f(g^{-1}x)=[\la(g)f](x)\quad \mbox{ for all}\  f \in
\H_\ell,\  g,x,\in G.
\]
It has been shown by Bernat-Pukanszky and Vergne that the unitary
representation $\pi_\ell$ is irreducible, and that $\pi_\ell$ is equivalent
to $\pi_{\ell'}$, if and only if $\ell $ and $\ell'$ lie on the same coadjoint
orbit, i.e. if and only if $\Ad^*(G)\ell=\Ad^*(G)\ell'$ (see \cite{Bernat} or
\cite [Theorem 8]{LeptinLudwig}) . Moreover, every irreducible unitary
representation of $G$ is
 equivalent to some $\pi_\ell$. This shows that one has a bijection
\[
K: {\g^*/\Ad^*(G)\to \hat G,\ \   \Ad^*(G)\ell\mapsto [\pi_\ell],}
\]
called the {\em Kirillov-map}. Here, $[\pi_\ell]$ denotes the equivalence class
 of $\pi_\ell$, and $\hat G$ the (unitary) dual of $G$, i.e. the set of all
equivalence classes of unitary irreducible representations of $G$.

\setcounter{equation}{0}
\section{The topology of $\hat G$}

Suppose again $G$ to be exponential, and denote by $C^*(G)$ the
{\em $C^*$-algebra of $G$}, which is, by definition, the completion of
$L^1(G)$
with respect to the {\em $C^*$-norm}
\[
||f||_{C^*}:=\sup\limits_{\pi\in \hat G} ||\pi(f)||, \quad f\in L^1(G).
\]
Since $G$ is amenable, $||f||_{C^*}$ is in fact equal to $||\la(f)|| $, where
$\la$ denotes the left-regular representation
(see (\cite{Pier}).

If $\pi\in \hat G$, $\pi$ extends uniquely to an irreducible unitary
representation
of $C^*(G)$, also denoted by $\pi$, and we let $I_\pi$ be the kernel of
$\pi$ in
 $C^*(G)$. This two-sided ideal is by definition a so-called
{\em primitive ideal}, and we denote by Prim$(G):=\{I_\pi:\pi\in \hat G\}$
the set of all primitive ideals of $C^*(G)$. We endow Prim$(G)$ with the
{\em Jacobson topology}. Thus a subset $C$ of Prim$(G)$ is closed if and
only if $C$ is the {\em hull} $h(I)$ of an ideal, i.e. $C=h(I)
:=\{J\in \Prim(G):J\supset I\}$. For any subset $A$ of $\Prim(G)$, we denote
 by $\ker A:=\bigcap _{J\in A} J$ the {\em kernel} of $A$,
which is an ideal in $C^*(G)$.

In any $C^*$-algebra, a closed two-sided ideal $I$ is always the kernel of its
hull, i.e.

\begin{equation}
I=\bigcap_{J\in \Prim C^*(G), J\supset I} J ;
\end{equation}
see e.g. \cite[2.9.7]{Dixmier}.

Now, since exponential Lie groups are so-called type I groups, the mapping
\[
\iota:\hat G \ni [\pi]\mapsto I_\pi\in \Prim(G)
\]
is a bijection (see \cite[\S 6]{LeptinLudwig}). In particular,
 $\iota\circ K:\g^*/\Ad^*(G)\to \Prim(G)$ is  bijective.

 Even more is true:
If we endow $\g^*/\Ad^*(G)$ with the quotient topology induced by the
topology  of
 $\g^*$, then $$\iota\circ K \mbox{ \em  is  a homeomorphism}$$ (see
\cite[\S 3,
Theorem 1]{LeptinLudwig}). We introduce on $\hat G$ a topology by pulling
back the
topology of  $\Prim(G)$  via $\iota$.

Our proof of Theorem 1 will make use of the following results, the first of
which is taken from \cite{Ludwig-Mueller}.

\begin{theorem}
Suppose $G$ is an exponential solvable Lie group, and let $\ell\in \g^*$. If
 the orbit $\Om(\ell)=\Ad^*(G)\ell$ is closed, then $\pi_\ell(C^*(G))$ is
the algebra
of  all compact operators on $\H_\ell$. In particular, $\pi_\ell(f)$ is
compact for
 every $f\in L^1(G)$.
\end{theorem}

The second result is a kind of ''Riemann-Lebesgue Lemma''. Let us call an
element
$\nu\in \g^*$ a {\em character}, if $\nu([\g,\g])=\{ 0\}.$

\begin{theorem}
Suppose $G$ is an exponential solvable Lie group, and  let $\ell\in \g^*$ with
coadjoint orbit $\Om:=\Om(\ell)$. Assume that the restriction of $\Om$
to the nilradical $\n$ of the Lie algebra $\g$ is closed. Then the orbit
$\Om$ is itself closed, and for any real character $\nu$ of $\g$  which
does not
vanish on the stabilizer $\g(\ell)$ of $\ell$, we have that
\begin{equation}\lim_{|\tau|\to\infty}\Om+\tau\nu=\lim_{|\tau|\to\infty}
\Om(\ell+\tau\nu)=\infty
\end{equation}
in the orbit space.
In particular,
\begin{equation} \lim_{|\tau|\to\infty} ||\pi_{\ell+\tau\nu}(f)||=0
\end{equation}
for every $f\in L^1(G).$
\end{theorem}

\bigskip

\noi{\bf Proof.}
Let $p:=\ell|_{\n}$ be the restriction of $\ell$ to $\n$.
The stabilizers $G(\ell)$ and $G(p)$ of
$\ell$ respectively of $p$ in $G$ are closed connected subgroups, and
we have $G(\ell)\subset G(p)$. There exists a closed subset $T$
of $G$ such that $G$
is the topological product of $T$ and $G(p)$, i.e. such that the mapping
$$T\times G(p)\to G,\ \ (t,u)\to t\cdot u,$$ is a homeomorphism.
In the same way let $S$ be a closed subset of $G(p)$ such that the mapping
$$S\times G(\ell)\to G(p),\ \ (s,u)\to s\cdot u,$$ is a homeomorphism
(see \cite {LeptinLudwig}). For $g\in G$ and $m\in \g^*$, let us write the
action of $g$
on $m$ as $$\Ad^*(g)m:=g\cdot m.$$
Put $\m(\ell):=\g(\ell)+\n.$ It is well-known that
\begin{equation}G(p)\cdot \ell=\ell+\m(\ell)^\perp.
\end{equation}

In fact, if $H=\exp{\h}:=G(p)$, then $\h=\{X\in\g: \ell([X,Y])=0 \quad\forall
Y\in\n\}.$ Therefore, if $Y\in\m(\ell)$ and $X\in\h$, then
$$\ell(e^{\ad X} Y)-\ell(Y)=p([X,Y]) +\frac 12 p([X,[X,Y]])+\dots =0,$$
since $\ad^*(X)p=0.$ This implies that $H\cdot\ell\subset \ell+\m(\ell)^\perp.$

Moreover, since the bilinear form $B_\ell(X,Y):=\ell([X,Y])$ is
non-degenerate on
$\g$ modulo $\g(\ell)$, we have
$$\g(p)/\g(\ell)\simeq (\n+\g(\ell))^\perp=\m(\ell)^\perp,$$
hence $\dim H\cdot \ell=\dim\m(\ell)^\perp.$ We thus obtain (2.4).

Assume now that
$\lim_{n\to\infty}\Om+\tau_n\nu=\Om(\ell')$, for some $\ell'\in \g^*$,
which means
that there exists a sequence $\{m_n\}_n=\{\ell_n+\tau_n \nu\}_n$ tending to
$\ell'$
in $\g^*$, where $\ell_n\in\Om$ and $\tau_n\in
\RR.$

We can write $\ell_n=(t_n s_n)\cdot \ell,$ with $t_n\in T$ and $s_n\in S.$
 Since $s_n\cdot \ell=\ell+q_n$ for some $q_n\in \m(\ell)^\perp$, it
follows that
$$\ell_n=(t_ns_n)\cdot \ell=t_n\cdot\ell+q_n,$$
and thus
$$(\ell_n+\tau_n \nu)|_{\n}=t_n\cdot \ell|_{\n},$$
hence
$$\ell'|_{\n}=\lim_{n\to \infty}t_n\cdot \ell|_{\n}.$$
 Since the restriction of $\Om$ to $\n$ is closed, we have that
$\ell'|_{\n}=t'\cdot p$ for some $t'\in T,$ and since $\Om|_{\n}$ is
homeomorphic to $G/G(p)\simeq T,$ it follows  that
$$\lim_{n\to\infty}t_n=t'.$$
Let us now take an element $U\in \g(\ell)$ such that $\nu(U)\ne 0.$
Then
$$\ell'(U)=\lim_{n\to\infty}t_n\cdot \ell(U)+\lim_{n\to\infty}\tau_n\nu(U)
=t'\cdot \ell(U)+\lim_{n\to\infty}\tau_n\nu(U).$$
Hence $\lim_{n\to\infty}\tau_n=\tau'$
exists, and it follows that the sequence $\{q_n\}_n$ convergences, hence also
$\lim_{n\to\infty}s_n=s'$ exists. Finally
$$\ell'=(t's')\cdot\ell+\tau'\nu\in \Om+\tau'\nu.$$
This proves (2.2), and (2.3) is an immediate consequence of (2.2) (see
\cite{Dixmier}).

\hfill Q.E.D.

\bigskip

\setcounter{equation}{0}
\section{Representations on mixed $L^p$-spaces}

We assume again that
$\g=\g_0\supset \g_1\supset \dots\supset \g_q=\n\supset \dots\supset
\g_m=\{0\}$ is a composition sequence passing through $\n$.
Let us assume that $\n$ is a nilpotent ideal containing $[\g,\g]$.

Let $\ell\in \g^*$, and let $\p(\ell)=\p$ be the  Vergne-polarization for
$\ell$
 associated to a fixed Jordan-H\"older sequence
\[
\g=\a_0\supset \a_1\supset \dots\supset \a_r=\{0\}
\]
refining this composition sequence. Then obviously $\p_0:=\p\cap \n$ is a
Vergne-polarization for $\ell_0:=\ell|_{\n}$.
As in the preceding proof, let $\g(\ell_0):=\{X\in\g: \ell_0([X,Y])=0 \quad
\forall Y\in\n\}$
be the stabilizer of $\ell_0$ in $\g$. Then
\begin{equation}
\p\subset \g(\ell_0) +\p_0.
\end{equation}
In fact, choose $k$ such that $\a_k=\n.$ Then, for $j\le k$ and
$X\in\a_j(\ell),$ we have
$\ell_0([X,Y])=0$ for every $Y\in\n,$ since $\n\subset\a_j$. This shows that
$\p\subset \g(\ell_0)+\sum_{j\ge k}\a_j(\ell_0)=\g(\ell_0)+\p_0.$


Next, for every $j\geq q$, we choose a subspace $\vv_j$ in $\g_j$ of
dimension $\le 2$,
such that $\g_j+\p_0=\vv_j\oplus(\g_{j+1} +\p_0) $, and define  the index set
$J$  as follows:
$$
J:=\{j\in\{q,\dots,m-1\}:\vv_j\ne \{0\}\}.
 $$
Write $J$  as an ordered  $d$-tuple
 $$
J=\{j_1 < \dots < j_d\},
$$
where $d:=\# J,$
and put $\ww_i:=\vv_{j_i}\subset \n, \ i=1,\dots,d$, and
$\ww:=\ww_1\oplus\dots
\oplus \ww_d$. We shall often identify $\ww$ with the direct product
$\ww_1\times\dots\times \ww_d$.

The space $\ww$ then forms a complementary subspace to the polarization
$\p_0$ in $\n,$
i.e.
\begin{equation}
\n=\ww\oplus\p_0.
\end{equation}

Let us  choose a linear subspace $\b$ of $\p$ such that
\begin{equation}
\p=\b\oplus \p_0.
\end{equation}
Then $\b\cap\n=\{0\}$, so that we may choose a subspace $\h$ of $\g$
containing $\n$
 such that
$$
\g:=\b\oplus\h.
$$
Then $\h$ is an ideal in $\g,$  and we may choose a subspace $\a$ of $\h$
such that
$$
\h=\a\oplus\n.
$$

Then we have $\p\cap\h=\p_0,$ and, by (3.3), (3.2),
\begin{equation}
\g=\a\oplus\b\oplus\n=\a\oplus(\p+\n)=\a\oplus\ww\oplus\p.
\end{equation}

Let $P:=\exp\p,$ $P_0:=\exp\p_0$ and $N:=\exp\n.$ Then the mapping
$$\Phi=\Phi_{G,P}:\a\times\ww\times P\to G,\quad
(S,(w_1,\dots,w_d),p)\mapsto \exp (S)\exp(w_1)\dots \exp(w_d)p,
$$
with $w_j\in\ww_j,$ is a diffeomorphism, and
$$
E=E_{G/P}: \a\times\ww\to G,\quad
(S,w)\mapsto \Phi(S,w,e)
$$
provides a {\em section} for $G/P$, i.e. $\a\times\ww\ni (S,w)\mapsto E(S,w)P$
is a diffeomorphism from  $\a\times\ww$ onto $G/P$.

Similarly,  $\ww\ni w\mapsto E(0,w)P_0$ is a diffeomorphisms from $\ww$
onto $N/P_0$.

We shall later make use of the ''global chart'' $E$ for $G/P$ in order to
construct
a more concrete realization of the induced representation $\pi_{\ell}$ on a
Euclidean $L^2-$
space, which will then also allow for the construction of more general
representations
on mixed $L^p$ spaces. Crucial for this construction will be the subsequent
analysis
of ''roots'' on $G$.

To begin with, let is  construct a  decomposition of $\p_0$ into subspaces
$\s_j$
 subordinate to our Jordan-H\"older sequence. To this end, choose  for
every $j\geq q$
 a subspace $\r_j$ in $\g_j$ of dimension $\le 2$,
such that $\g_j\cap\p_0=\g_{j+1}\cap\p_0 \oplus \r_j$, and  define another
index set
$I$ as follows:
$$
I:=\{j\in\{q,\dots,m-1\}:\r_j\ne \{0\}\}. $$
Write $I$  again as an ordered $e$-tuple $$
I=\{j'_1 < \dots < j'_e\},
$$
where $e:=\# I,$
and put $\s_i:=\r_{j'_i}, \ i=1,\cdots,e$ . Then
$$\p_0=\s_1\oplus\dots\oplus\s_e\simeq \s_1\times\dots\times\s_e,$$
and the mapping
$$\Phi_P:\b\times\p_0\to P,\quad (T,Y_1,\cdots,Y_e)\mapsto
\exp(T)\exp (Y_1)\cdots \exp(Y_e)\in P$$
is a diffeomorphism which identifies the Lebesgue measure on
 $\b\times\p_0$ with the Haar measure on $P$.

Define also for every $j=q,\cdots,m-1$ the subspace $\uu_j$ of $\g_j$ by
$\uu_j:=\r_j+\vv_j.$ Then $\uu_j$ is the direct sum

$$\uu_j=\r_j\oplus\vv_j,$$
and
\begin{equation}\ \g_j=\uu_j\oplus\g_{j+1},\quad j=q,\dots,m-1.
\end{equation}
In particular, we have
$$\n=\uu_q\oplus\cdots\oplus\uu_{m-1}.$$
According to (3.5),  for $j=q,\dots,m-1,$  $X\in \g $ and $U\in \g_j,$ we may
write
\begin{equation}
\ad(X)(U)=\al_j(X)U+U_j,
\end{equation}
where $U_j$ is the component of $\ad(X)(U)$ in $\g_{j+1},$ and where
$\al_j(X)$ is an
endomorphism of $\uu_j$. Then $\al_j$ is an irreducible representation of
$\g$ on
 $\uu_j,$ which we shall call a {\em root} of $\g$.
Since $G$ is exponential, the eigenvalues of $\al_j(X),$ considered as an
endomorphism
 of the complexification of $\uu_j,$ are of the form $\al(1+i\beta)$, where
$\al$
and $\beta$ are real numbers.
For $X\in \g$ and $j=q,\dots,m-1,$ let
$$\tau_j(X):=\tr \ \ad_{\g_j/\g_{j+1}}(X)=\tr\  \al_j(X),$$
where by $\ad_{\g_j/\g_{j+1}}(X)$ we denote the factorized adjoint action
of $X$ on
the quotient space $\g_j/\g_{j+1}.$

The functionals $\tau_j$ are characters of $\g$, since
$\ad_{\g_j/\g_{j+1}}(X)=0$
for every $X\in \n$. Since $\ad(\p)$ acts on $\g/\p,$ one finds that for
$X\in\p$ the
corresponding ''trace of $\ad(X)$ modulo $\p$'' is given by

\begin{equation}
\tr\  \ad_{\g_j+\p/\g_{j+1}+\p}(X)= \ve_j \tau_j(X), \quad j=q,\dots,m-1,
\end{equation}
where
$$\ve_j:=\frac{\dim \vv_j}{\dim \uu_j}, \quad j=q,\dots,m-1.$$

Observe that $\ve_j\ne 0$ if and only if $\g_j+\p/{\g_{j+1}+\p}\simeq
\vv_j$ is non-trivial,
 i.e. if and only if $j\in J=\{j_1 < \dots < j_d\}.$ For $i=1,\dots,d$ and
$T\in\p$
 we shall therefore put $\la_i(T):=\al_{j_i}(T),$ so that for every
$w_i\in\ww_i$

\begin{equation}
\ad(T)(w_i)=\la_i(T)w_i\ \hbox{ modulo }\ \g_{j_i+1},\quad i=1,\dots,d.
\end{equation}
Then, by (3.7) and (3.8), we have
\begin{equation}
\tr\ \ad_{\g/\p}(T)=\sum_{i=1}^d \ve_{j_i}\tr\,\la_i(T), \quad T\in\p.
\end{equation}

Observe now that also the mapping
$$\Psi:\b\times\a\times N\to G, \quad (T,S,n)\mapsto \exp(T)\exp(S) n,$$
is a diffeomorphism. And, for  every  $R\in\a,$ $w=(w_1,\dots,w_d)\in \ww$ and
 $S\in \a $, $T\in\b\subset\p$, $n\in\n$  we have
\begin{eqnarray}
&(&\exp(T) \exp(S) n)^{-1}E(R,w)\nonumber\\
&=&n^{-1}\exp (-S)\exp (T)^{-1}\exp (R)\exp (T)
\left(\prod_{i=1}^d \exp(e^{-\ad (T)} w_i)\right) \,\exp(-T).
\end{eqnarray}
{}From (3.8) and (3.10), one can  deduce that
\begin{equation}
(\exp(T) \exp(S) n)^{-1}E(R,w)=E(R-S,\om(R,w,T,S,n)) p(R,w,T,S,n)^{-1},
\end{equation}
where $\om:\a\times\ww\times \b\times \a\times N\to \ww,\quad
p:\a\times\ww\times \b\times \a\times N\to P$ are analytic
mappings which depend polynomially on $w$ and $n$, and where
$\om=(\om_1,\dots,\om_d)$,
 with

\begin{equation}
\om_i(R,w,T,S,n)=e^{-\la_i(T)} (w_i)+\tilde \om_i(R,w_1,\dots,w_{i-1},T,S,n).
\end{equation}
Because of (3.3) and (3.10), we have  $p(R,w,T,S,n)=\exp (T)\ \hbox{mod}\
P_0$, i.e.
\begin{equation}
p(R,w,T,S,n)=\exp(T)\,\nu(R,w,T,S,n),
\end{equation}
with $\nu(R,w,T,S,n)\in P_0=P\cap N$.

Putting $p:=\exp (T)\in P,$ we therefore obtain
\begin{eqnarray*}
\De_{G,P}&(p(R,w,T,S,n))=\De_{G,P}(p)=
\frac{ [\det e^{\ad_\g(T)}]^{-1}} {[\det e^{\ad_\p(T)}]^{-1}}\\
=&e^{-\tr\, \ad_{\g}(T)+\tr \, \ad_{\p}(T)} =e^{-\tr \, \ad_{\g/\p}(T)},
\end{eqnarray*}
hence, by (3.9),
\begin{equation}
\Delta_{G,P}(p(R,w,T,S,n))=
e^{-\sum\limits_{i=1}^d \ve_{j_i}\tr\, \la_i(T)} .
\end{equation}

In particular, if we define the real character $\De$ on $G$ by
$$\De(\exp(X)):=\exp\left( -\sum_{j=q}^{m-1}\ve_j \tau_j(X)\right)=
\exp\left( -\sum_{i=1}^d\ve_{j_i} \tau_{j_i}(X)\right), \quad X\in \g,
$$
then
\begin{equation}
\De_{G,P}(p)=\De(p) \ \hbox{ for every } \ p\in P.
\end{equation}

Now, in oder to realize the representation on a Euclidean $L^2$-space, we
first
observe that the left-invariant linear functional $\oint_{G/P} \,d\dot x$
in (1.1) is
given by
\begin{equation}%
\oint_{G/P} f(x)\,d\dot x =\int_{\a\times\ww} f\circ E(R,w)\,dR\, dw \quad
\forall f\in \E(G,P)
\end{equation}
(see \cite[Theorem 2]{LeptinLudwig}).
For $f\in\E (G,P),$ let us put
$$\tilde f(x):=\De(x) f(x), \quad x\in G.$$
Since, by (3.15), $\De$ is a character of $G$ which extends $\De_{G,P}$
from $P$ to $G$,
 we have for $x\in G$ and $p\in P$
$$\tilde f(xp)=\De(xp)f(xp)=\De(x)\De(p)\De_{G,P}(p)^{-1}f(x)=\tilde f(x).$$
Thus the mapping $\A^2:f\mapsto \tilde f$ is a linear isomorphism from
$\E(G,P)$
onto the space
\begin{eqnarray*}\tilde\E (G,P):=\{f\in C(G,\CC): f\
 \mbox{ has compact support modulo}\ P,\\
 \ \mbox {and}\  f(xp)=f(x)\quad \forall x\in G, p\in P\}.\\ \end{eqnarray*}
Identifying functions on $G/P$ with $P$-right-invariant functions on $G$,
we thus see
that $\tilde\E(G,P)\simeq C_0(G/P).$  Moreover, if we define for any
continuous function
$f$ with compact support on $G/P$ its integral by

$$\int_{G/P}f(x)d\dot x:=\int_{\a\times\ww}\De^{-1}(E(R,w)) f(E(R,w))\,
dR\, dw,$$
then clearly
\begin{equation}
\oint_{G/P} f(x)\,d\dot x =\int_{G/P} \tilde f\,d\dot x \quad \hbox{for every}\
f\in\E(G,P).
\end{equation}
Comparing with (1.1), we find in particular that
\begin{equation}
\int_G F(x)\,dx =\int_{G/P}\int_P F(xp)\De(xp)\, dp\, dx \quad \hbox{for
every}\
f\in C_0(G).
\end{equation}

Let us define the space
\begin{eqnarray*}\tilde\E(G,P,\ell):=\{f\in C(G,\CC): f \ \mbox{ has
compact support
 modulo}\  P, \\
\ \mbox{and} \ f(xp)=\o{\chi_\ell(p)}  f(x)
\ \forall x\in G,& p\in P\},\end{eqnarray*}

\noindent  endowed with the norm $||\cdot||_2$ given by
\begin{equation}
||f||_2^2:=  \int_{G/P} |f(x)|^2\,d\dot x=
\int_{\a\times\ww}|\De^{-1/2}(E(R,w)) f(E(R,w))|^2\, dR\, dw.
\end{equation}
Observe that $|f|^2\in \tilde \E(G,P),$ if $f\in \tilde\E(G,P,\ell).$ Let
$\tilde\H_\ell$ denote
the completion of $\tilde\E(G,P,\ell)$ with respect to this norm. It is
obvious that
the mapping  $\A:f\mapsto \De^{1/2}f$ is a linear isomorphism between
$\E(G,P,\ell)$ and $\tilde \E(G,P,\ell),$ which extends to an isometric
isomorphism
of the Hilbert space $\H_\ell$ onto the Hilbert space $\tilde\H_\ell.$
Moreover,
$\tilde{\H}_\ell$ is nothing but the space

\begin{eqnarray*}L^2(G/P,\ell):=\{f:G\to\CC: f\  \hbox{is measurable,}
\ \hbox{and} \ f(xp)=\o{\chi_\ell(p)}  f(x)\ \\
\hbox{for a.e.}\ x\in G \ \hbox{and every}\ p\in P,\
\hbox{s.t.}\  ||f||_2<\infty\}.
\end{eqnarray*}
We may therefore intertwine the representation $\pi_\ell$ with the operator
$\A$
 in order to obtain a unitarily equivalent representation $\tilde\pi_\ell$ on
$L^2(G/P,\ell)$, given by
$$\tilde\pi_\ell(g):=\A\pi_\ell(g)\A^{-1}, \quad g\in G.$$
A straight-forward computation shows that $\tilde\pi_\ell$ is given
explicitly by
$$
[\tilde\pi_\ell(g)f](x)=\De(g)^{1/2}f(g^{-1}x) \ \hbox{for all}\ f\in
L^2(G/P,\ell),
\ g,x\in G,
$$
i.e.
\begin{equation}
\tilde\pi_\ell(g)=\De(g)^{1/2}\la(g), \quad g\in G.
\end{equation}

For a ''multi-exponent'' $\u{p}=(p,p_1,\dots,p_d)\in [1,\infty[^{1+d}$, let us
now define the
{\it mixed $L^{p}$-space}\  $L^{\u p}(G/P,\ell)\simeq
L^p(\a,(L^{p_1}(\ww_1,L^{p_2}(\ww_2,\dots)))$ by

\begin{eqnarray*}L^{\u p}(G/P,\ell):=\{f:G\to\CC: f\  \hbox{is measurable,}
\ \hbox{and} \ f(xp)=\o{\chi_\ell(p)}  f(x)\ \\
\hbox{for a.e.}\ x\in G \ \hbox{and every}\ p\in P,\
\hbox{s.t.}\  ||f||_{\u p}<\infty\},
\end{eqnarray*}
where the mixed $L^{\u p}$-norm is given by
$$
||f||_{\u p} $$
$$:=
\left(\int_{\a}
\left( \int_{\ww_1} \dots
\left( \int_{\ww_{d-1}}
\left( \int_{\ww_d}
\left|(\De^{-\frac 1p} f)(E(R,w_1,\dots,w_d))
\right|^{p_d}dw_d
\right)^{\frac{p_{d-1}}{p_d}} dw_{d-1}
\right)^{\frac{p_{d-2}}{p_{d-1}}}\dots dw_1 \right)^{\frac p{p_1}} dR
\right)^{\frac 1p}.
$$
The space $\tilde\E(G,P,\ell)$ is dense in $\Lp,$ for any $\u p$.

Put
$$\ga_i(\u p):=\frac 1{p_i}, \quad i=1,\dots, d,$$
and define the character $\de_{\u p}$ of $\g$ by
$$\de_{\u p}(X):=\sum\limits_{i=1}^{d}
\ga_i(\u p)\ve_{j_i}\tau_{j_i}(X),\quad \hbox{if}\ X\in \p,$$
 $$\de_{\u p}(X):=\frac 1p \ \sum\limits_{j=q}^{m-1} \ve_j\tau_j(X), \quad
\hbox{if}\ X\in \a+\n.$$
Observe that $\de_{\u p}$ is well-defined, since $(\a+\n)\cap\p\subset\n,$
and since
$\tau_j$ vanishes on $\n.$ The corresponding character of $G$ is given by
$$
\De_{\u p}(\exp (X)):=e^{-\de_{\u p}(X)},\quad X\in\g.$$

Notice that for $\o 2:=(2,\dots,2)$ we have
 \begin{equation}\De_{\o 2}=\De^{1/2}\ \hbox{ and } \ ||\cdot||_{\o 2}=
||\cdot||_2. \end{equation}
Observe also that for $R\in \a, w\in\ww,$ we have
$$\De(E(R,w))=\De(\exp(R)).$$
Let  $T\in\b$ and  $R\in\a, w\in\ww$. Then,  by (3.11), (3.12),  we have
$$\exp(T)^{-1} E(R,w)=E(R,\{e^{-\la_i(T)} w_i+\tilde \om_i(R,w_1,
\dots,w_{i-1},T,0,e)\}_{i=1}^d)
\, p(R,w,T,0,e)^{-1}.$$
{}From the definition of the norm $||\cdot||_{\u p}$ we therefore obtain
 $$||\la(y)f||_{\u p}=\De_{\u p}(y)^{-1}||f||_{\u p}\quad \hbox{for every}
 \ f\in\Lp,\,y\in P.$$
Similarly, for $S\in\a$, we have that
$$\exp S^{-1} E(R,w)=E(R-S,\{
w_i+\tilde\om_i(R,w_1,\dots,w_{i-1},0,e)\}_{i=1}^d)
p(R,w,0,e)^{-1},$$
hence
$$||\la(\exp S)f||_{\u p}=\Delta_{\u p}(\exp (S))^{-1} ||f||_{\u p} \quad
\hbox{for every} \ f\in\Lp,\,S\in \a.$$
 Finally, if we choose $n\in N$, then of course
$$||\la(n)f||_{\u p}=||f||_{\u p}.$$
It is now clear that  we obtain  an {\em isometric} representation
$\pi_{\ell}^{\u p}$
of $G$ on $\Lp$ by letting
$$[\pi_{\ell}^{\u p}(g)f](x):=\De_{\u p}(g) f(g^{-1}x),\quad g,x\in
G,f\in\Lp,$$
i.e.
\begin{equation}\pi_{\ell}^{\u p}(g)=\De_{\u p}(g) \la(g),\quad g\in G.
\end{equation}

Notice that by (3.20) and (3.21), the representation $\pi_{\ell}^{\o 2}$ is
unitarily
 equivalent to $\pi_\ell,$ i.e.
\begin{equation}
\pi_{\ell}^{\o 2}\simeq \pi_\ell.
\end{equation}
In the sequel, we shall work with $\pi_{\ell}^{\o 2}$ in place of $\pi_\ell.$
With a slight abuse of notation, we shall therefore denote $\pi_{\ell}^{\o 2}$
simply by $\pi_\ell.$  Observe that then for every function $f\in L^1(G)$ such
that $\De_{\u p}\De^{-1/2} f\in L^1(G),$ the
operator $\pi_{\ell}^{\u p}(f)$ is given by the formula
\begin{equation}
\pi_{\ell}^{\u p}(f)=\pi_\ell(\De_{\u p}\De^{-1/2} f),
\end{equation}
acting boundedly on the space $\Lp.$ More generally, we have

\begin{prop}
Let
${\u p},{\u q}\in [ 1,\infty [ ^{1+d},$ and let $f\in L^1(G)$ such that
$\De_{\u p}\De_{\u q}^{-1} f\in L^1(G).$ Then the operator $\pi_{\ell}^{\u
p}(f)$
extends uniquely from $\Lq\cap\Lp$ to a bounded operator on $\Lq,$ given by
the
formula
\begin{equation}
\pi_{\ell}^{\u p}(f)=\pi_\ell^{\u q}(\De_{\u p}\De_{\u q}^{-1} f).
\end{equation}
In particular, one has
\begin{equation}
||\pi_{\ell}^{\u p}(f)||_{\Lq\to\Lq}\le ||\De_{\u p}\De_{\u q}^{-1} f||_1.
\end{equation}
\end{prop}


\section{Auxiliary results}

\setcounter{equation}{0}
\subsection{Compact operators acting on mixed  $L^p$-spaces  }

In this subsection, let $M:=X\times Y$ be a product of two measure spaces
$(X,dx)$ and
$(Y,dy)$. For $1\le p<\infty, $ denote by $L^{\u p}$ the mixed $L^p$-space
$L^p(X,L^2(Y)), $ endowed with the norm
$$||f||_{\u p}:=
\left(\int_{X}(\int_Y| f(x,y)|^2\, dy)^{p/2}\, dx\right)^{1/p}.$$ By
$\Id_A$ we denote
the indicator function of a set $A.$

\setcounter{subsection}{0}
\begin{lemma}
(i) Let $E=\{E_1,\dots, E_n\}$ be a family of disjoint measurable subsets
of finite
measure in  $X,$  and denote by $S: L^{\u p}\to L^{\u p}$ the associated
averaging
operator
$$ S(f)(x,y):=\sum_{j=1}^n \frac 1 {|E_j|}\left( \int_{E_j}f(u,y)\, du\right)
\, \Id_{E_j}(x),\quad f\in L^{\u p},\  (x,y) \in X\times Y$$
with respect to the first variable. Then the operator norm of $T$ is
bounded by 1,
for every $p.$

(ii) Similarly, let  $\F=\{F_1,\dots, F_m\}$ be a family of disjoint
measurable subsets of
 finite measure in $Y,$  and
denote by $T: L^{\u p}\to L^{\u p}$ the associated  averaging operator
$$  T(f)(x,y):=\sum_{j=1}^m \frac 1 {|F_j|}\left(
\int_{F_j}f(x,v)\, dv\right)\, \Id_{F_j}(y),\quad f\in L^{\u p},\  (x,y)
\in X\times Y$$
with respect to the second variable.
Then the operator norm of $T$ is bounded by 1, for every $p.$
\end{lemma}

\bigskip

\noi{\bf Proof.} In order to prove (i),  we write $f_x(y):=f(x,y).$ Then
$$S(f)_x=\sum_j\frac 1{|E_j|} \Id_{E_j}(x) f_j,$$
where  $f_j:=\int_{E_j}f_u\, du\in L^2(Y).$ Therefore
$$||S(f)_x||_2=\sum_j\frac 1{|E_j|}||f_j||_2\Id_{E_j}(x), $$
hence
$$\int_X ||S(f)_x||_2^p\, dx = \sum_j\frac 1{|E_j|^p}||f_j||_2^p |E_j|
=\sum_j |E_j|^{1-p}||f_j||_2^p.$$
But, by Minkowski's integral inequality and H\"older's inequality,
\begin{eqnarray*}
||f_j||_2^p&=&||\int_{E_j} f_u\, du||_2^p\le \left( \int_{E_j}||f_u||_2\,
du\right)^p\\
&\le& |E_j|^{p/p'}\int_{E_j}||f_u||_2^p\, du=
|E_j|^{p-1}\int_{E_j}||f_u||_2^p\, du.
\end{eqnarray*}
Consequently,
$$||S(f)||_{\u p}^p\le \sum_j\int_{E_j}||f_u||_2^p\, du=\int_X||f_u||_2^p\, du
=||f||_{\u p}^p.$$

The proof of (ii) is even simpler. Indeed, for $f\in L^{\u p}$, we have
\begin{eqnarray*}
||T(f)||_{\u p}^p&=&\int_X\left(\int_Y\left| \sum_j \frac 1 {|F_j|}
\Id_{F_j}(y)\int_{F_j}
f(x,v)\, dv\right| ^2\, dy\right)^{p/2}dx\\
&=&\int_X\left(\sum_j\frac 1 {|F_j|} \left| \int_{F_j}f(x,v)\, dv \right|^2
\right)^{p/2}
dx \le \int_X\left(\sum_j{\int_{F_j}|f(x,v)|^2\, dv}\right)^{p/2}dx\\
&\le& \int_X\left(
\int_{Y}|f(x,v)|^2\, dv\right)^{p/2}\, dx=||f||_{\u p}^p.
\end{eqnarray*}

\bigskip
\hfill Q.E.D.

\medskip
We are now in a position to prove the following variant  for mixed
$L^p$-spaces of an
interpolation theorem by Krasnoselskii \cite{Krasnoselskii}.

\begin{theorem}
Let $p,q\in[1,\infty[, \ p\ne q,$ and let $K$ be a linear operator on
$L^{\u p}+L^{\u q}$
which maps $L^{\u p}$ compactly into $L^{\u p}$ and $L^{\u q}$ boundedly
into $L^{\u q}.$
Then $K$ is a compact operator from  $L^{\u r}$ to $L^{\u r},$ for every
$r$ lying
strictly between $p$ and $q$.
\end{theorem}

\bigskip

\noi{\bf Proof.} Observe that the space $\cS$ of simple functions of the form
$\sum_j\al_j\Id_{E_j}\otimes\Id_{F_j}$, where the $E_j$ and $F_j$ form
finite families
of measurable subsets of finite measure in $X,$ respectively $Y,$ lies
dense in
$L^{\u r},$ for every $1\le r<\infty.$

Denote by $B_1(0)$ the unit ball centered at the origin in $L^{\u p}.$ Since
 $\K:=\o{K(B_1(0))}$ is a compact subset of $L^{\u p},$ for every $k\in\NN$
 we may thus find simple functions $f_1^k,\dots,f_{j_k}^k$ in $\cS$ such
that,
for any $g\in\K,$ there exists a $j$ such that $||g-f_j^k||_{\u p}<\frac
1{k+1}.$

Choose next finite families $\E_k=\{E_1^k,\dots,E_{n_k}^k\}$, respectively
$\F_k=\{F_1^k,\dots,F_{m_k}^k\}$ of disjoint measurable subsets of $X$,
respectively of
$Y$, such that every $f=f_j^k$  can be written as a linear combination of
functions
 of the form $\Id_{E_r^k}\otimes \Id_{F_s^k}.$

 Denote by $S_k $ and  $T_k$ the averaging operator associated to $\E_k ,$
respectively
  $\F_k$ in Lemma 4.1,
 and let $R_k$ be the operator of finite rank given by
$$R_k:=S_k\circ T_k,\quad k\in \NN.$$
Then $R_k f_j^k=f_j^k$ for any $j$, and so, if  $g\in \K,$ and if we choose
$j$ such that
$||g-f_j^k||_{\u p}<\frac 1{k+1},$ then we obtain from Lemma 4.1 that
$$||g-R_kg||_{\u p}\le ||g-f_j^k||_{\u p}+||R_k(f_j^k-g)||_{\u p}\le \frac
2{k+1}.$$
This shows that $\lim_{k\to \infty}R_k(g)=g$ in $ L^{\u p}$ for every $g\in
\K$. As a
consequence, we obtain
\begin{equation}
\lim_{k\to\infty}||R_k\circ K-K||_{L^{\u p} \to L^{\u p}}=0.
\end{equation}

For, otherwise we could find a sequence of functions $f_k$ in $B_1(0)$ and
$\ve>0,$  such
that
\begin{equation} ||f_k||_{\u p}=1 \ \hbox{ and }\ ||(R_k-\Id)\circ K
f_k||_{\u p}
\ge \ve \quad  \forall k.
\end{equation}
Passing to a subsequence, if necessary, we could then assume that the
sequence of functions
$g_k:=Kf_k$ had a limit $g$ in $\K,$ since $\K$ is compact. This would imply
$\lim_{k\to\infty}(R_k-\Id)(g_k-g)=0,$ hence
$\lim_{k\to\infty}(R_k-\Id)\circ K(f_k)=0,$ contradicting (4.2).

On the other hand
\begin{equation}
||R_k\circ K-K||_{L^{\u q}\to L^{\u q}}\leq 2||K||_{L^{\u q}\ \to L^{\u
q}},\quad \forall
k\in \NN,
\end{equation}
 because $||R_k||_{L^q\to L^q}$ is bounded by 1 for every $k$. Applying
the  Riesz-Thorin interpolation theorem, it follows from (4.1) and (4.2) that
$$\lim_{k\to\infty}||R_k\circ K-K||_{L^{\u r}\ \to L^{\u r}}=0,$$
for every $r$ lying strictly between $p$ and $q$. Our assertion follows.

\bigskip
\hfill Q.E.D.
\medskip

In the sequel, we shall denote the space of compact operators on a Banach
space $E$
by $\K(E).$ Moreover, if $K$ is an operator as in Thm. 4.2, and if we
consider $K$
as a compact operator from $L^{\u r}$ to $L^{\u r},$ for $r$ lying  strictly
between $p$ and $q,$ then we shall also write $K_r$ in place of $K$. The
spectrum
of $K_r\in\K(L^{\u r})$ will be denoted by $\si_r(K).$ Given $p\in
[1,\infty[,$
we shall denote the conjugate exponent by $p',$ i.e. $\frac 1p+\frac 1{p'}=1.$

\medskip
As for the  spectra of $K$ on different $L^{\u r}$-spaces, we  have

\begin{prop} 
Let $1<p_0\le 2,$  and let $K$ be a linear operator on
$L^{\u p_0}+L^{\u p_0'},$ mapping the space $L^{\u p_0},$ as well as the
dual space
$L^{\u p_0'},$ compactly into itself, so that, by the preceding theorem,
$K$ is a
compact operator on $L^{\u p}$ for every $p\in [p_o,p_o'].$ Assume further
that
$K$ is self-adjoint on $L^{\u 2}$.
 Let $\la\in\si_{p_o}(K)\setminus \{0\}.$ Then $\la$ is real, and every
generalized
 eigenvector of $K_{p_o}$ associated to $\la$ is in fact an eigenvector,
lying in
$\bigcap_{p\in [p_0,p_0']}L^{\u p}.$ In particular, all $L^{\u p}$- spectra
coincide, i.e.
$\si_p(K)=\si_2(K),$ and so do the  eigenspaces corresponding to non-zero
eigenvalues,
 for every $p\in[p_0,p_0']$.
\end{prop}

\bigskip

\noi{\bf Proof.} Let $\la\in\si_{p_o}(K)\setminus \{0\}, $ and
let $E\subset D$ be two compact neighbourhoods of the line segment
$[\la,{\o \la}]$, which
are invariant under complex conjugation, such that $D$ is also a
neighbourhood of $E$.
Since the non-zero eigenvalues are isolated, by shrinking $E$ and $D$, if
necessary,
we may assume that $\si_{q}(K)\cap D\subset [\la,{\o\la}], $
for $q=p_0,2,p_0'$. Then for any $\mu\in D\setminus E^\circ$, we have that
$(K-\mu)^{-1}$ exists on $L^{\u q}$ for $q=p_0,2,p_0',$ and  there exists a
constant
$C>0$, such that
$$||(K-\mu)^{-1}||_{L^{\u p_0}\to L^{\u p_0}} +||(K-\o\mu)^{-1}||_{L^{\u p_0'}
\to L^{\u p_0'}}|| \le C,\quad \forall \mu\in D\setminus E^\circ.$$
Thus, by interpolation, we obtain
\begin{equation}
||(K-\mu)^{-1}||_{L^{\u p}\to L^{\u p}} \le C,\ \hbox{for every}\  \mu\in
D\setminus E^
\circ, \ p\in[p_0,p_0'].
\end{equation}
In particular $(K_p-\mu)^{-1}$ exists on $L^{\u p}$.
This implies that no point $\mu\in D\setminus [\la,\o\la] $ lies in  any of
the sets
$\si_p(K), \ p\in[p_0,p_0'],$  i.e. $\si_p(K)\cap D \subset [\la,\o \la]$. Let
$$\D:=\bigcap_{p\in [p_0,p_0']} L^{\u p}.$$
Then $\D$ is dense in $L^{\u p}$ and invariant under $K_p,$ for every $p\in
[p_0,p_0'].$

 We show that, for $\mu\in D\setminus [\la,\o\la],$ the restriction
$$(K_p-\mu)^{-1}|_{\D}\ $$ of $(K_p-\mu)^{-1}$ to the joint core $\D$ does
not depend on
 $p\in[p_0,p_0'].$

Indeed, let $\xi\in L^{\u p}\cap L^{\u q},$ and let
$$\eta_p:=(K_p-\mu)^{-1}\xi,\quad \eta_q:=(K_q-\mu)^{-1}\xi,$$
for given  $p,q\in [p_0,p_0']$. Then, for $\vp\in \D$, we have that
$$\l \eta_p,(K-\o\mu)\vp\ra =\l (K_p-\mu)^{-1} \xi, (K-\o\mu)\vp \ra=\l
\xi,\vp \ra
=\l \eta_q, (K-\o\mu)\vp\ra.$$
Since $(K-\mu)$ is invertible on $L^{\u p}$, it follows that $(K-\o\mu)$ is
invertible on the dual space $(L^{\u p})',$ and the same applies to $(L^{\u
q})'$.
Consequently,
$(K-\o\mu)(\D)$ is dense in $(L^{\u p})'$ as well as in $(L^{\u q})',$ so
that $\eta_p$
and $\eta_q$ coincide as linear functionals on $(L^{\u p})'$ and on $(L^{\u
q})',$
and so $\eta_p=\eta_q\in L^{\u p}\cap L^{\u q}$.
This implies that $(K_q-\mu)^{-1}\xi=(K_{p_0}-\mu)^{-1}\xi$ for every
$\xi\in \D$ and
every $q\in[p_0,p_0'].$

Let $\Gamma$ be a simple curve in $D\setminus E^\circ$ whose winding number
with respect
to every point in the interval $[\la,\o\la]$ is one,  and let
$$P_p:=\int_{\Gamma} (K_p-\ga)^{-1}\, d\ga,\quad  p\in[p_0,p_0'].$$
Then the  operators $P_p$ coincide on $\D$ and  are bounded by a common
constant,
by (4.4). The image  $\im P_p$ of $P_p$ is the sum of the generalized
eigenspaces of $K_p$ over
all eigenvalues of $K_p$ which are contained in $\si_p\cap [\la,\ \o\la]$.
In particular, the image of $P_p$ is finite dimensional,  since $K_p$ is
compact. Since
$\im P_p=\o{P_p(\D)},$ we thus have
 $$\im P_p=P_p ( \bigcap_{q\in [p_0,p_0']} L^{\u q})=\im P_2.$$
In particular, $\im P_p$ does not depend on $p$.

 Since $K$ is self-adjoint on $L^{\u 2}$, its restriction to $\im P_2$ is
also self-adjoint
 on $\im P_2,$ and so every eigenvalue of $K,$ when acting  on $\im P_p$,
is real, for any $p\in[p_0,p_0'].$ This applies  in particular to $\la. $
We thus see
 that $\si_p(K)\cap D=\{\la\}$,
that $ \im P_2\subset \D,$ and that $\im P_2$ is in fact the eigenspace of
$K$ associated to $\la,$ for any $p\in[p_0,p_0'].$
Furthermore, every eigenvector of $K$ for the eigenvalue $\la$ in $L^{\u
p}$ is
contained in $\bigcap_{q\in [p_0,p_0']}L^{\u q}.$

\bigskip
\hfill Q.E.D.

\setcounter{subsection}{1}
\subsection{Approximate units of Herz-Schur multipliers }

Let $G$ be a locally compact  group, endowed with a left-invariant
Haar measure $dx$. If $K$ is a continuous function on $G$, we shall
 denote by $\la (K):C_0(G) \to C(G)$ the convolution operator  given by
$$\la(K)(\vp):=K\star \vp,\quad \vp\in C_0(G);$$
here, $\la$ denotes again the left-regular representation.
In case that $\la(K)$ extends to a bounded operator on $L^p(G),$ for some
$p\in[1,\infty[,$
 we shall denote the (unique) extension also by $\la(K).$

The following result is a  consequence of the well-known theory of
Herz-Schur multipliers
and Fourier-Figa-Talamanca-Herz-Eymard  algebras, see e.g. \cite{Herz},
\cite{Eymard},
\cite{Figa}. For the convenience of the reader, we shall  provide a proof.

\setcounter{subsection}{3}
\begin{lemma} Let $K\in C(G),$ and assume that $\la(K)$ is  bounded
on $L^p(G).$ Let  $\psi\in L^{p'}(G)$ and $ \eta\in L^{p}(G).$ Then also
$\la((\psi\star\check \eta)K)$ is bounded on $L^p(G),$ and
\begin{equation}
||\la((\psi\star\check\eta)K)||_{L^p\to L^p}\le
||\psi||_{p'}||\eta||_{p}||\la(K)||_{L^p\to L^p}.
\end{equation}
\end{lemma}

\bigskip

\noi{\bf Proof.}
If $f,g\in C_0(G),$ then
\begin{eqnarray*}
\l \la((\psi\star\check\eta)K) f,g\ra &=& \int\int K(y)
(\psi\star\check\eta)(y)f(y\inv x) g(x)\, dy\,
dx\\
&=& \int\int\int K(y)\psi(t)\check\eta(t\inv y)f(y\inv x) g(x)\,dy\, dx\, dt\\
&=&\int\int\int K(y)\psi(xt)\eta(y\inv xt)f(y\inv x) g(x)\,dy\, dx\, dt\\
&=& \int \l K\star (\eta_t f),\psi_t g\ra\, dt,
\end{eqnarray*}
where we have used the abbreviation  $h_t(x):=h(xt).$ Thus, by H\"older's
inequality and
Fubini's theorem, we obtain
\begin{eqnarray*}
&&|\l \la((\psi\star\check\eta)K) f,g\ra |\\
&\le&\int ||\la(K)||_{L^p\to L^p}||\eta_t f||_p
||\psi_t g||_{p'}\, dt\\
&\le& ||\la(K)||_{L^p\to L^p}\left( \int\int |\eta(yt)|^p |f(y)|^p\, dy\,
dt\right)^{1/p}
\left( \int\int |\psi(xt)|^{p'} |g(x)|^{p'}\, dx\, dt\right)^{1/p'}\\
&=&||\la(K)||_{L^p\to L^p}||\eta||_p ||\psi||_{p'} ||f||_p ||g||_{p'}.
\end{eqnarray*}
This implies (4.5).

\bigskip
\hfill Q.E.D.
\medskip

Assume next that $G$ is a  Lie group.  If $T$ is a right-invariant bounded
operator on
$L^p(G),$ then it follows from the Schwartz' kernel theorem that there
exists a unique
distribution $K\in\D'(G),$ such that
\begin{equation}
T\vp=K\star \vp=\la(K)(\vp),\quad \hbox{for every}\ \vp\in\D(G).
\end{equation}

We shall then also denote $K$ by $T\de_e.$  Inversely, given a distribution
$K\in\D'(G),$
we shall denote by $\la(K):\D(G) \to C^\infty(G)$ the convolution operator
given by
$$\la(K)(\vp):=K\star \vp,\quad \vp\in\D(G).$$
In case that $\la(K)$ extends to a bounded operator on $L^p(G),$ for some
$p\in[1,\infty[,$
 we shall again  denote also the (unique) extension by $\la(K).$ The following
proposition makes use of  ideas in \cite{Leptin}, \cite{Eymard}.

\begin{prop}
Let $G$ be  an amenable Lie group, let $1\le p<\infty,$ and let
$K\in\D'(G),$ such that
the convolution operator $\la(K)$ is bounded on $L^p(G).$ Assume further
that $\la(K)\in
C^*(G).$ Then $\la(K)$ is bounded on $L^q(G)$ for every exponent $q$ lying
between
$2$ and $p$, and there exists a sequence $\{K_n\}_n$ in $C_0(G),$ such that
\begin{equation}
\lim_{n\to\infty} ||\la(K)-\la(K_n)||_{L^2\to L^2}=0,
\end{equation}
and
\begin{equation}
||\la(K_n)||_{L^q\to L^q}\le ||\la(K)||_{L^q\to L^q}\quad \forall n\in\NN,
\end{equation}
for every $q$ lying  between $2$ and $p$.
\end{prop}

\bigskip

\noi{\bf Proof.}
Assume first that $K\in C(G).$
We choose an increasing sequence $\{A_n\}_{n\in \NN}$ of compact subsets of $G$
such that
$G=\bigcup_{n\in \NN}A_n$. Let $\ve>0.$ Since $G$ is amenable, we can find
compact
subsets $U_n$ of $G$ such that
$$\frac {|(A_n\inv U_n)\bigtriangleup U_n|} {|U_n|}<\ve.$$
Put
$$\phi_{n,\ve}:= \frac 1 {|U_n|} \Id_{U_n}\star\check\Id_{U_n}\in C_0(G).$$
For any $q\in[1,\infty[,$ let us put
$$\psi_{n,\ve}=\frac 1 {|U_n|^{1/p'}} \Id_{U_n}, \quad
 \eta_{n,\ve}:=\frac 1 {|U_n|^{1/p}} \Id_{U_n},\quad n\in\NN.$$
 Then of course
$$||\psi_{n,\ve}||_{p'}=1=||\eta_{n,\ve}||_{p}, \quad n\in\NN,$$
and
$$ \phi_{n,\ve}=\psi_{n,\ve}\star \check\eta_{n,\ve}.$$
If $a\in A_n,$ then
\begin{eqnarray*}
||\la(a) \eta_{n,\ve}-\eta_{n,\ve}||_{p}^{p} &=&
\frac 1 {| U_n|}\int_G \left| \Id_{U_n}(a\inv t)-\Id_{U_n}(t)\right|^{p}\, dt\\
 &=& \frac 1 {|U_n|}\int_G \Id_{(aU_n)\bigtriangleup U_n}(t)\, dt \le \ve.
\end{eqnarray*}

Observe that $\phi_{n,\ve}(e)=1.$ Therefore,  for any $a\in A_n,$
 we have that
\begin{eqnarray*}
|\phi_{n,\ve}(a)-1|&=&|\phi_{n,\ve}(a)-\phi_{n,\ve}(e)|
=|\int \psi_{n,\ve}(t)[\eta_{n,\ve}(a\inv t)-\eta_{n,\ve}]\, dt|\\
&\le& ||\psi_{n,\ve}||_{p'} ||\la(a)\eta_{n,\ve}-\eta_{n,\ve}||_{p}\le
\ve^{1/p}.
\end{eqnarray*}

Thus, if we put $\phi_n:=\phi_{n,1/n},\ n\ge 1,$ then we see that the
functions $\phi_n$
tend to $1$, uniformly on compacta. Let us put $K_n:=\phi_n K.$ Then
$K_n\in C_0(G),$
and, by the preceding proposition, the $K_n$ satisfy (4.8). Moreover, since
$\la(K)\in
C^*(G),$ there exists a sequence $\{f_j\}_j$ in $C_0(G)$ such that
$\lim_{j\to\infty} ||\la(K)-\la(f_j)||=0,$ where by $||\cdot||$ we denote
the operator
norm on $L^2(G).$ Since $\lim_{n\to\infty} ||f_j-\phi_n f_j||_1=0,$ for
every $f_j$
we  see that $\lim_{n\to\infty} ||\la(f_j)-\la(\phi_n f_j)||=0.$ Moreover,
by Lemma 4.4,
 we have
\begin{eqnarray*}
||\la(K)-\la(K_n)||&\le& ||\la(K)-\la(f_j)||+||\la(f_j)-\la(\phi_n f_j)||+
  ||\la(\phi_n(f_j-K))||\\
&\le 2& ||\la(K)-\la(f_j)||+||\la(f_j)-\la(\phi_n f_j)||.
\end{eqnarray*}
Thus, given $\ve>0,$ we may first choose $j$ such that
$||\la(K)-\la(f_j)||<\ve/2, $
 and then $n_0\in\NN^\times$ such that $||\la(f_j)-\la(\phi_n f_j)||<\ve/4$
for every $n\ge n_0.$ Then
$||\la(K)-\la(K_n)||<\ve$ for every $n\ge n_0,$ i.e. (4.7) is also satisfied.

For an arbitrary $K\in\D'(G)$ satisfying the assumptions of Prop. 4.5, we
may argue as follows.
   We fix an approximate identity $\{\chi_i\}_i$ in $\D(G)$  such that
$||\chi_i||_1=1$
for every $i,$ and put $K^{(i)}:=\chi_i\star K.$ Then $K^{(i)}\in C(G).$
Moreover,
\begin{eqnarray*}
||\la(K)-\la(K^{(i)})||&\le& ||\la(K)-\la(f_j)||+||\la(f_j)-\la(\chi_i\star
f_j)||+
  ||\la(\chi_i\star(f_j-K))||\\
&\le 2& ||\la(K)-\la(f_j)||+||f_j-\chi_i\star f_j||_1,
\end{eqnarray*}
where $\{f_j\}_j$ is again a sequence in $C_0(G)$ such that
$\lim_{j\to\infty} ||\la(K)-\la(f_j)||=0.$
Thus, arguing similarly as before, we see that
$$\lim_{i\to \infty}||\la(K)-\la(K^{(i)})||=0.$$
Let us put $K^{(i)}_\nu:=\phi_\nu K^{(i)}\in C_0(G).$ Then, we know already
that
$$||\la(K^{(i)}_\nu)||_{L^q\to L^q}\le ||\la(K^{(i)})||_{L^q\to L^q}
\le ||\chi_i||_1 ||\la(K)||_{L^q\to L^q}=||\la(K)||_{L^q\to L^q}.$$
Moreover, we have
$$||\la(K)-\la( K^{(i)}_\nu)||\le ||\la(K)-\la( K^{(i)})||+||\la(K^{(i)})
-\la( K^{(i)}_\nu)||.$$
Thus, if $\ve>0$ is given, we may choose $i$ such that $||\la(K)-\la(
K^{(i)})||<\ve /2,$
and subsequently $n_0\in\NN^\times$ such that $||\la(K^{(i)})-\la(
K^{(i)}_\nu)||<\ve/2$ for
every $\nu\ge n_0.$ Then $||\la(K)-\la( K^{(i)}_\nu)||< \ve$ for every
$\nu\ge n_0.$ This
shows that we may find a subsequence $K_n$ among the $K^{(i)}_\nu$
satisfying (4.7) and
(4.8).

\bigskip
\hfill Q.E.D.

\setcounter{equation}{0}
\section{An analytic family of compact operators}

Let us now choose a fixed  sub-Laplacian $L$ on $G,$ and denote by
$\{e^{-tL}\}_{t>0}$ the heat semigroup generated by $L.$ We recall some
well-known facts
about this semi-group (see e.g. \cite{Ludwig-Mueller}, also for further
references).

 For every $t>0,$  $e^{-tL}$ is a convolution operator
\begin{equation}
e^{-tL} f=h_t\star f,
\end{equation}
where the $\{h_t\}_{t >0}$ form a 1-parameter semigroup of smooth
probability measures in $L^1(G).$

Moreover, as a consequence of Gaussian estimates for
the heat kernels, one has the following extension of Lemma 5.1 in
\cite{Ludwig-Mueller}, whose proof carries over to the present situation.

\begin{prop}
Let $\s$ be a subspace of $\g$ complementary to $\n,$ for instance
$\s=\a+\b,$ so that
the mapping $\s\times N\ni (S,n)\mapsto \exp(S) n\in G$ is a diffeomorphism
from $\s\times N$
onto $G,$ and fix a norm $|\cdot |$ on $\s.$
For any $a\ge 0, j\in \NN$, put
\[
h^{a,j}_1(\exp(S)n):=|S|^je^{a|S|} h_1(\exp(S)n),\hskip1cm (S,n)\in \s\times N.
\]
Then $h_{a,j}\in L^1(G)$. Moreover, there is a constant $C_a > 0$,  such that

\begin{equation}
||h^{a,j}_1||_1 \le C_a^{j+1} \Gamma \left( \frac j 2 +1\right).
\end{equation}
\end{prop}
\bigskip

If $\chi$ is any continuous, real or complex character of $G$, with
differential $d\chi\in\g_\CC^*,$
then $\chi(\exp(S) n)=e^{d\chi(S)}$. We therefore have the following

\begin{cor} $\chi h_1\in L^1(G)$ for every continuous character $\chi$ of $G$.
\end{cor}

{}From now on, we shall make the following
\medskip

\noi{\bf Assumption.} $\ell\in\g^*$ satisfies Boidol's condition, and
$\Om(\ell)|_\n$
is closed.

\noi Moreover, we assume   and  $p\in [1,\infty[,\ p\ne 2,$ is fixed.

\medskip

Then, since $\ell$ satisfies (B), there exists at least one root $\la$ of
$\g$, such that
 $\la|_{\g(\ell)}$ is a non-trivial  root of the $\g(\ell)$-module $\g/\p$
(see \cite{Boidol}). Consequently, there exists at least one index $i\in
\{1,\dots,d\},$
such that  $$\la_i|_{\g(\ell)}\ne 0,$$
 and $(\g_{j_i}+\p)/(\g_{j_i+1}+\p)\ne \{ 0 \}.$ Notice that the latter
condition is
equivalent to $\ve_{j_i}\ne 0.$
Choose $i$ minimal with these properties, and put

\begin{equation}
\o p:=(p,\dots,p,2,2,\dots, 2),
\end{equation}
where the last $p$ occurs at the $i$-th position.

 Then, for $T\in\p,$ we have
$$(\de_{\o p}-\de_{\o 2})(T)=\sum_{k=1}^{i}\frac 1p \ve_{j_k} \tau_{j_k}(T)
+\sum_{k=i+1}^{d}\frac 1 2 \ve_{j_k} \tau_{j_k}(T)- \sum_{k=1}^{d}\frac 1 2
\ve_{j_k}
 \tau_{j_k}(T)
=(\frac 1p -\frac 12)\sum_{k=1}^i\ve_{j_k}\tau_{j_k}(T),$$

and for $X\in \a+\n$ one has
$$(\de_{\o p}-\de_{\o 2})(X)=(\frac 1 p-\frac 1 2)\sum_{k=1}^{m-1} \ve_k
\tau_k(X).$$

Thus, if we define the real character $\nu$ of $\g$ by
$$\nu (X):=\sum_{k=1}^i\ve_{j_k}\tau_{j_k}(X),\quad \hbox{if}\ X\in \p,$$
 $$\nu (X):=\sum_{k=1}^{m-1} \ve_k \tau_k(X), \quad
\hbox{if}\ X\in \a+\n,$$
then
\begin{equation}
\De_{\o p}\De_{\o 2}^{-1}(\exp(X))=e^{(\frac 12 -\frac 1p)\nu(X)}, \quad
X\in \g.
\end{equation}
Moreover, since $\tau_{j_k}(T)=0$ for $1\le k<i$ and $T\in\g(\ell),$ we have
\begin{equation}
\nu|_{\g(\ell)}=\ve_{j_i} \tau_{j_i}|_{\g(\ell)}\ne 0.
\end{equation}

For any complex number $z$ in the strip
$$\Si:=\{\ze\in\CC: |\Im\ze|<1/2\},$$
let $\De_z$ be the complex character of $G$ given by
$$\De_z(\exp(X)):=e^{-i z\nu(X)}, \quad X\in\g,$$
and $\chi_z$ the unitary character
$$\chi_z(\exp(X)):=e^{-i\Re (z)\nu(X)}, \quad X\in\g.$$
Since, by (5.4),
$$\De_z(x)=\chi_z \De_{\o {p(z)}}\De_{\o 2}^{-1},$$
if we define $p(z)\in ]1,\infty[$ by the equation
\begin{equation}
\Im (z)=1/2 -1/p(z),
\end{equation}
 we see that the representation $\pi_\ell^z$, given by
\begin{equation}
\pi_\ell^z(x):=\De_z(x)\pi_\ell(x)=\chi_z(x)\pi_\ell^{{\o {p(z)}}}(x),
\quad x\in G,
\end{equation}
is an isometric representation on the space $L^{{\o {p(z)}}}(G/P,\ell).$

Observe that for $\tau\in\RR,$ we have $p(\tau)=2,$ and $\pi_\ell^\tau=
\chi_\tau\otimes\pi_\ell$ is a unitary representation on $L^{\o 2}$. Moreover,
\begin{equation}\pi_\ell^{\tau}\simeq\pi_{\ell-\tau\nu,}
\end{equation}
since the mapping  $f\mapsto {\o \chi}_\tau f$ intertwines the representations
$\chi_\tau \otimes\pi_\ell$ and $\pi_{\ell-\tau\nu}.$

Let us put
$$T(z):=\pi_\ell^z(h_1)=\pi_\ell(\De_z h_1), \quad z\in\Si,$$
and let us shortly write
$$L^{\o p}:={L^{\o p}(G/P,\ell)}, \quad 1\le p<\infty,$$
where $\o p$ is given by (5.3).

Since, by (3.24),
\begin{equation}
T(z)=\pi_\ell^{\o q}(\De_z \De_{\o 2}\De_{\o q}^{-1}h_1),
\end{equation}
it follows from Cor.5.2 and Prop. 3.1 that the operator $T(z)$ leaves
$L^{\o q}$
invariant for every $1\le q<\infty,$ and is bounded on all these spaces.
Much more is even true.  Let us write $T_q(z)$ in place of $T(z),$
 if we consider $T(z)$ as a bounded operator on $L^{\o q}.$ The spectrum of
$T_q(z)$
will be denoted by $\si_q(z).$

\begin{prop}
For every $q\in ]1,\infty[,$ the mapping $\Si\ni z\mapsto T_q(z)$ is an
analytic
family of compact operators in the sense of Kato (\cite{Kato}). Moreover,
if $\tau\in\RR,$
then $T_2(\tau)$ is self-adjoint on $L^{\o 2} $.
\end{prop}
\bigskip

\noi{\bf Proof.}  By Thm. 2.2, the orbit $\Om(\ell)$ is closed. Moreover,
Cor. 5.2 shows
that $\De_z h_1\in L^1(G), $ and consequently $T(z)=\pi_\ell(\De_z h_1)$ is
a compact
operator on $L^{\o 2}.$ On  the other hand, in a similar way we see from
(5.9) that
$T(z)$ is bounded on $L^{\o q},$ for every $1\le q<\infty.$  Since $L^{\o
q}$ is a mixed $L^p$-space
of the type $L^p(X,L^2(Y))$ considered in Section 4.1, we may apply Thm.
4.2 to
conclude that $T_q(z)$ is compact for every $q\in ]1,\infty[.$

Next, let $\ze\in\Si$ be fixed, and consider $T(\ze +z),$ for $|z|$
sufficiently small.
We have
 \begin{equation}
T(\ze+z)=\pi_\ell(\De_{\ze+z}h_1)=\sum_{j=0}^\infty\frac{(-iz)^j}{j!}S_j,
\end{equation}
where
$$S_j:= \pi_\ell( (\nu\circ \log)^j\De_\ze h_1)=
\pi_\ell^{\o q}((\nu\circ \log)^j\De_\ze \De_{\o 2}\De_{\o q}^{-1}h_1).$$
By Prop. 5.1 we see that
 \begin{equation}
||S_j||_{L^{\o q}\to L^{\o q}}\le C_q^{j+1}\Gamma(\frac j2+1),
\end{equation}
 where the constant $C_q>0$ stays bounded whenever $q$ runs through a compact
interval. Moreover, arguing for $S_j$ as we did for $T(z)$ before, one
finds that $S_j\in
\K (L^{\o q}),$ for every $q\in ]1,\infty[.$ Thus, by (5.10) and (5.11),
the mapping
$z\mapsto T_q(z)$ is holomorphic from $\Si$ into $\K (L^{\o q}),$ for
$1<q<\infty $ (it
even extends to an entire mapping from $\CC$ into $\K (L^{\o q}).$)

Finally, if $\tau\in\RR,$ then $\pi_\ell^\tau$ is a unitary representation
on $L^{\o 2},$
so that
$T(\tau)^*=\pi_\ell^\tau(h_1)^*=\pi_\ell^\tau(h_1^*)=\pi_\ell^\tau(h_1)=T(\tau).
$

\bigskip
\hfill Q.E.D.

We can now prove the following perturbation result.

\begin{prop} Let $1\le p_0<2.$ There exist an open neighborhood $U$ of a point
$z_0\in\RR$ in the complex strip $\Si$ and holomorphic mappings
$$\la:U\to\CC$$
and
$$\xi:U\to \bigcap_{p_0\le p\le p_0'} L^{\o p},$$
such that $\xi(z)\ne 0$ and
\begin{equation}
T(z)\xi(z)=\la(z)\xi(z)\quad\hbox{for every}\quad z\in U.
\end{equation}
Moreover, shrinking $U,$ if necessary, one can find a constant $C>0$ such that
\begin{equation}
||\xi(z)||_{L^{\o p}}\le C\quad \hbox{for every }\quad z\in U,\ p\in[p_0,p_0'].
\end{equation}
\end{prop}

\bigskip

\noi{\bf Proof.} Let $z_0\in\Si$ be real. Then, by Prop. 5.3 and Prop 4.3,
the $L^{\o p}$-
spectrum of $T_p(z_0)$ is independent of $p,$ and agrees thus with
$\si_2(z_0).$ Moreover,
for every non-trivial eigenvalue $\la_0\in\si_2(z_0)$ of $T(z_0),$  every
generalized
eigenvector $\xi_0$ associated to $\la_0$ is in fact an eigenvector, lying in
$\D:=\bigcap_{p_0\le p\le p_0'}L^{\o p}.$ Let us for instance choose for
$\la_0$ the largest
eigenvalue of $T(z_0).$

 Consider the three analytic families
$$\{T_q(z)\}_{z\in\Si}, \quad\hbox{for}\quad q=p_0,2, p_0'.$$
By choosing $z_0$ in such a way that $\la_0$ is a non-branching eigenvalue
for all
three analytic families and applying analytic perturbation theory (see
\cite{Kato}), we
may find an open connected neighborhood
$U$ of some real point $z_0\in\RR$ in $\Si$ and three holomorphic families
$U\ni z\mapsto \la_q(z)$ of eigenvalues for the operator $T_q(z),\quad
q=p_0,2, p_0',$
which all coincide at $z_0$ with the eigenvalue $\la_0.$ Moreover, we may
choose a
neighborhood $W$ of $\la_0$ such that
$$W\cap \si_q(z)=\{ \la_q(z)\}\quad\hbox{ for}\quad  z\in U,\ q=p_0,2,p_0'.$$

 Let
$$P_q(z):=\int_\Gamma (T_q(z)-\mu)\inv \, d\mu, \quad z\in U,$$
where $\Gamma$ is a circle in $W$ winding around $\la_0$ once.
By shrinking $U$, if necessary, we may also assume that the curve $\Gamma$
separates
$\la_q(z)$ from the remaining elements of $\si_q(z),$ for every $z\in U$ and
$q=p_0,2, p_0'.$  Then $P_q(z)$ projects onto the generalized eigenspace
$E_q(z)$ of
$T_q(z)$ associated with the eigenvalue $\la_q(z)$, for $q=p_0,2,p_0'.$
Moreover, $\D$ is the  core considered in the proof
of Prop.4.3, and  we had seen there that, for real $z\in U,$ the
restrictions of the
resolvents $(T_p(z)-\mu)\inv$ to $\D$ do not depend on $p,$ so that the
same applies to
the projectors $P_p(z).$ Since $P_q(z)$ depends holomorphically on $z$, it
follows that
\begin{equation}
P_q(z)|_\D=P_2(z)|_\D\quad\hbox{for every}\quad z\in U\cap\RR,\  q=p_0,2, p_0'
\end{equation}
Moreover, since we may assume that $P_q(z)$ is uniformly bounded on $L^{\o
q},$ for
$z\in U$ and $q=p_0,2, p_0',$ by interpolation we derive from (5.14) that
\begin{equation}
||P_2(z)\xi||_{L^{\o p}}\le C||\xi||_{L^{\o p}} \quad \forall z\in U,\
p\in[p_0,p_0'],
\ \xi\in\D.
\end{equation}
Choose now $\xi_0\in\D\setminus\{0\} $ such that $T(z_0)\xi_0=\la_0\xi_0,$
and put
$$\xi(z):=P_2(z)\xi_0, \quad z\in U.$$
Then $\xi(z)\in\D\setminus\{ 0\}$ for every $z\in U,$ (5.13) holds because
of (5.15),  and
the mapping $z\mapsto \xi(z)\in\D$ is
holomorphic with respect to every $L^{\o p}$-norm, $p_0\le p\le p_0',$
provided we choose
$U$ sufficiently small.

Furthermore, for real $z\in U,$ we have $T(z)\xi(z)=\la(z)\xi(z), $ and
since both sides of
this equation depend holomorphically on $z,$ it remains valid for every
$z\in U.$

\bigskip
\hfill  Q.E.D.


\setcounter{equation}{0}
\section{The proof of Theorem 1}

We have to show that, under the assumptions made in the previous section,
there
exist  a point $\la_0$ in the $L^2$-spectrum of $L$ and an open
neighborhood $\U$ of
$\la_0$ in $\CC,$ such that every $L^p$-multiplier $F\in C_\infty(\RR)$
extends
holomorphically to $\U$. Since $\si_2(L)=[0,\infty [=:\RR_+,$ we may
restrict ourselves to
multipliers on $\RR_+.$

\begin{lemma}
Let $1\le p<\infty,$ and let $F\in \M_p(L)\cap C_\infty(\RR_+).$ Then
$F(L)\in C^*(G),$
and $F(L)$ is bounded on $L^q(G)$ for every $q$ lying between $p$ and $p'.$
Moreover, there
exists a constant $C\ge 0,$ such that
\begin{equation}
||F(L)||_{L^q\to L^q}\le C, \quad \hbox{if}\quad |\frac 1q-\frac 12|\le
|\frac 1p-\frac 12|.
\end{equation}
\end{lemma}

\bigskip

\noi{\bf Proof.}
Since  $F\in C_\infty(\RR_+),$ there exists a sequence of functions of
Laplace transform
type $\tilde \vp_n(\la)=\int_0^\infty \vp_n(t)e^{-\la t}\, dt,$ where
$ \vp_n \in C_0(\RR_+),$ which converges uniformly to $F$ (see
\cite{Ludwig-Mueller},
Prop.2.1). Moreover, if $\vp$ is a real valued function on $\RR_+,$ then
$\tilde \vp$ is
real valued too, as is the convolution kernel
$\tilde \vp (L)\de_e=\int_0^\infty\vp(t)h_t\, dt$  associated to
$\tilde\vp(L).$
Let $\al_n:=\Re \vp_n,$ $\beta_n:=\Im \vp_n, $ and $F_1:=\Re F, F_2:=\Im F.$
Then $F_1$ is the uniform limit of the $\tilde\al_n,$ and $F_2$ is the
uniform limit of the
 $\tilde\beta_n$,
so that consequently, for every real-valued function  $f\in \D(G),$ one has
$$F_1(L)f=\lim_{n\to\infty}\tilde\al_n(L)f\ \hbox{and}\
F_2(L)f=\lim_{n\to\infty}\tilde\beta_n(L)f $$
in $L^2(G).$
This shows that the $F_1(L)f$ and $F_2(L)f$ are real-valued functions, whence
$$||F(L)f||_p=\left(\int_G (|[F_1(L)f](x)|^2+|[F_2(L)f](x)|^2)^{p/2}\,
dx\right)^{1/p}\ge
\max (||F_1(L)f||_p, ||F_2(L)f||_p).$$
Hence $F_1$ and $F_2$ are $L^p$-multipliers for $L$ too, and since
$F_1(L)^*=\o F_1(L)=F_1(L)$ as well as $F_2(L)^*=F_2(L),$ we see that $F_1(L)$
and $F_2(L)$ are also  bounded on $L^{p'}(G)$. This shows that
$F(L)=F_1(L)+iF_2(L)$ is
$L^{p'}$-bounded, and (6.1) follows by interpolation.

\bigskip
\hfill Q.E.D.

\bigskip

In view of Lemma 6.1, we may and shall assume in the sequel that $1<p<2.$
Let $K:=F(L)\de_e$
be the convolution kernel of $F(L),$ so that $F(L)\vp=K\star
\vp=\la(K)\vp,$ for
$\vp\in\D(G).$
According to Prop. 4.5,  choose a sequence $\{K_n\}_n$ in $C_0(G)$ such that
\begin{equation}
\lim_{n\to\infty} ||\la(K)-\la(K_n)||_{L^2\to L^2}=0,
\end{equation}
and
\begin{equation}
||\la(K_n)||_{L^q\to L^q}\le ||\la(K)||_{L^q\to L^q}\quad \forall n\in\NN,
\end{equation}
for every $q$ lying between $2$ and $p$.

We  now apply the transference theorem in order to
conclude
that
\begin{equation}
||\pi_\ell^z(K_n)||_{L^{\o{p(z)}}\to L^{\o{p(z)}}}\le ||\la(K_n)||
_{L^{p(z)}(G)\to L^{p(z)}(G)} \quad \forall n\in\NN, z\in\Si.
\end{equation}

To this end, recall that our representation $\pi_\ell^z $ acts
isometrically on a mixed
$L^p$-space of the type $L^p(X,L^2(Y)),$ with $p=p(z).$ Such a space can be
embedded into an $L^p$-space. Namely, $L^2(Y)$ is isometrically isomorphic
to a subspace  of $L^p(Z)$ (see \cite{MarcZygmund} Lemme 1 or \cite{Herzt}
Corollary 1). Then  of course $L^p(X,L^2(Y))$ is isometrically isomorphic to
a subspace of
$L^p(X\times Z)$.  However, the proof of Theorem 2.4  in
\cite{CoifmanWeiss} remains
valid also for bounded representations on closed subspaces of $L^p$-spaces,
and we can
thus apply the transference theorem in \cite{CoifmanWeiss} to obtain (6.4).

\medskip
\noi {\bf Remark.} (6.4) is an immediate consequence of the  inclusion
$A(\xi)A_p\subset A_p$ of Theorem A  in \cite{Herzt} and Theorem 6
in \cite{Herza}. It seems that Herz understood transference very well, but
in his publications rather concentrated on  abstract results (notably  in
\cite{Herzt}) and  made explicit only few consequences. We felt
that in case of (6.4)  (as in section 4) re-proving certain more or less
known   results is more convenient   then explaining how to translate
and properly combine  known results to get what we need.

\medskip

{}From (6.1), (6.3) and (6.4) we get
\begin{equation}
||\pi_\ell^z(K_n)||_{L^{\o{p(z)}}\to L^{\o{p(z)}}}\le C \quad \hbox{for every}
\quad z\in\Si_p,
\end{equation}
where $\Si_p$ denotes the smaller strip
$$\Si_p:=\{\ze\in\CC:|\Im \ze|<\frac 1p-\frac 12\}.$$

Next, letting $p_0:=p,$ choose holomorphic families $\la(z)$ of eigenvalues for
 $\pi_\ell^z=T(z)$
and associated eigenfunctions $\xi(z), \ z\in U,$  as in Prop. 5.4, and
assume w.r.
that $U\subset \Si_p.$

Fix $\psi\in C_0^\infty (G/P,\ell)$ such that $\l \xi(z), \psi\ra\ne 0$ for
every $z\in U$ (if
necessary, we have to shrink $U$ another time to achieve this),
and consider the complex functions
$$h_n:U\to \CC, \quad z\mapsto \l \pi^z_\ell(K_n)\xi(z), \psi\ra,\quad
n\in\NN,$$
which are holomorphic in $U.$

By (6.5) and (5.13) we obtain
$$|h_n(z)|\le ||\pi_\ell^z(K_n)||_{L^{\o{p(z)}}\to L^{\o{p(z)}}}\,
||\xi(z)||_{\o{p(z)}}
\, ||\psi||_{\o{p(z)'}}\le C,\quad \forall z\in U, n\in\NN.$$
The family of functions $\{h_n\}_n$ is thus a normal family of holomorphic
functions, and
so, by the theorems of  Montel and Weierstra\ss, there exists a subsequence
$\{h_{n_j}\}_j$ which converges locally uniformly to a holomorphic limit
function $h$ on $U.$

Now, if $z\in U$ is real, then the representation $\pi^z_\ell$ is unitary,
and since, by
(6.2), $\la(K_n)$ converges to $\la(K)$ in $C^*(G),$ as $n$ tends to $\infty$,
it follows that $ h_n(z)$ converges  on $U\cap \RR$ to  $\l
\pi_\ell^z(K)\xi(z),\psi\ra,$ i.e.
\begin{equation}
h(z)=\l \pi_\ell^z(K)\xi(z),\psi\ra,\quad z\in U\cap\RR.
\end{equation}
Let $\mu(z):=-\log \la(z), \ z\in U,$ where $\log$ denotes the principal
branch of the
logarithm. Then $\mu$ is holomorphic on $U,$ and from
$\pi_\ell^z(h_1)\xi(z)=T(z)\xi(z)
=\la(z)\xi(z)$ we obtain that $d\pi_\ell^z(L)\xi(z)=\mu(z)\xi(z).$ From
(6.6) we therefore
get (compare \cite{Ludwig-Mueller})
$$h(z)=\l \pi_\ell^z(F(L))\xi(z),\psi\ra =\l F(d\pi_\ell^z(L))\xi(z),\psi\ra
=F(\mu(z))\, \l \xi(z),\psi\ra,$$
i.e.
\begin{equation}
F\circ\mu(z)=\frac{h(z)}{\l \xi(z),\psi\ra}, \quad \forall z\in U\cap\RR.
\end{equation}

Now, clearly the right-hand side of (6.7) extends holomorphically to $U,$
and thus
$F\circ \mu$ extends to a holomorphic function on $U.$

However, from Thm. 2.2 and formula (5.8) we deduce that
$$\lim_{|\tau|\to\infty}||\pi_\ell^\tau(h_1)||=0,$$
hence $\la(\tau)\to 0$ and $\mu(\tau)\to +\infty$ as $\tau\to\infty.$ Thus
the function
$\mu$ is not constant, and so, varying the point $z_0\in U\cap\RR$
slightly, if necessary,
we may assume that $\mu'(z_0)\ne 0.$ Then $\mu$ is a local bi-holomorphism
near $z_0,$
and thus $F$ has a holomorphic extension to a complex neighborhood of
$\mu(z_0).$

\bigskip
\hfill Q.E.D.


\setcounter{equation}{0}
\setcounter{section}{6}
\section{An example of a closed orbit whose restriction to the nilradical is
non-closed}

Let $\g$ be the Lie algebra spanned by the basis
$$\B=\{R,S,T,X,Y,Z,M_1,M_2,N_1,N_2,\},$$
with non-trivial brackets given by
\begin{eqnarray*}
&&[T,X] = -X,\ [T,Y]=Y,\ [X,Y]=Z,\\
&& [R,T] =  M_1,\ [R,M_1]=M_2,\ [R,M_2]=-M_2, \\
&& [S,T] = N_1,\ [S,N_1]=N_2,\ [S,N_2]=N_2.
\end{eqnarray*}
Let $\al,\ \beta\in \RR\setminus \{0\},$ and denote by $\ell$ the element of
 $\g^*$ for which
$$\ell(Z)=1,\ \ell(M_2)=\al,\ \ell(N_2)=\beta, \quad\ell(U)=0$$
for all other elements\  $U$ of the basis $\B.$

The stabilizer of $\ell$ in $\g$ is the subspace
$$\g(\ell)=\span \{T,\ Z,\ N_1-N_2,\ M_1+M_2\}.$$
Let
$$g(x,y,r,s,n,m):=\exp{x X}\exp{yY}\exp {r R} \exp{s S}\exp{m M_2}\exp{n
N_2}.$$
Then the coadjoint orbit $\Om$ of $\ell$ is the subset
$$\Om=\{\Ad^*(g(x,y,r,s,m,n)\inv)\ell: \ x,y,r,s,m,n\in\RR\}.$$
Denote by
$$\B^*:=\{R^*,S^*,T^*,X^*,Y^*,Z^*,M_1^*,M_2^*,N_1^*,N_2^*,\}$$
the dual basis of $\B$. Then
\begin{eqnarray*}
\Om&=&\{-mR^*+nS^*+(\al(e^{-r}-1+r)+\beta(e^s-1-s)-xy)T^* +
\al(1-e^{-r})M_1^*\\
 &+&\al e^{-r}M_2^*+\beta(e^s-1)N_1^* +\beta e^s N_2^* -yX^*+x Y^* +Z^*:\
x,y,r,s,m,n\in\RR\}.
\end{eqnarray*}
We see that the restriction of $\Om$ to the nilradical
$\n=\span\{X,Y,Z,M_1,M_2,N_1,N_2\}$ is not closed, since, letting $r$ tend
to $+\infty$
and $s$  to $-\infty,$  the other parameters being fixed, one finds that
the functionals
$$\al M_1^*+\beta N_1^*-yX^*+x Y^*+Z^*$$
lie  in the closure of $\Om|_\n$.

If $\al/\beta>0$, then $\al(e^{-r}-1+r)+\beta(e^s-1-s)-xy$ tends to
infinity, provided
  $r$ tends to $+\infty$ and $s$ tends to $-\infty$ whereas $x,y$ stay
bounded. Hence the
orbit $\Om$ is closed, whenever $\al/\beta >0$.

On the other hand,  $\Om$ is not closed if $\al/\beta<0$, since then the
element
$\al M_1^*+\beta N_1^*$ is contained in the closure of the orbit ( take
$s:=\frac{\al+\beta+\al r} \beta $ and let $r$ tend to $+\infty$.)

It is easy to see that the subspace $\p:=\span \{T,Z,Y,M_1,M_2, N_1,N_2\}$
is the Vergne
polarization for $\ell$ associated to the composition sequence
\[
\g=\g_0\supset \g_1\supset\dots\supset \g_{10}=\{0\},
\]
where $\g_j$ is spanned by the $j+1$-st to last element or the ordered basis
$$R,S,T,X,Y,Z,M_1,M_2,N_1,N_2$$
of $\g.$

The root $\la$ associated to the quotient space
$$\span \{X,Y,M_1,M_2,N_1,N_2,Z\}/\span \{Y,M_1,M_2,N_1,N_2,Z\}$$
is the linear functional $\nu:=-T^*$. In particular $\ell$  does not satisfy
Boidol's condition (B).

 To fix the ideas, assume now for instance that $\al>0,\ \beta>0$ (which
means that $\Om$
 is closed), and consider a real sequence $\{\tau_k\}_{k\in \NN}$ such that
$\lim_{k\to \infty}\tau_k=+\infty$. Then the sequence of orbits
$\Om_k:=\Om+\tau_kT^*,\ k\in \NN,$ tends to infinity in the orbit space.
Indeed, otherwise,
for a  subsequence, also indexed by $k$ for simplicity of notation, for
every $k$ there
would exist an element
\begin{eqnarray*}
\ell_k&=&(-m_kR^*+n_kS^*+ \al(e^{-r_k}-1+r_k)+\beta(e^{s_k}-1-s_k)-x_k
y_k+\tau_k)T^* \\
&+&\al(1-e^{-r_k})M_1^*
+\al e^{-r_k}M_2^*+\beta (e^{s_k}-1)N_1^* +\beta e^{s_k} N_2^* -y_kX^*+x_k
Y^* +Z^*\ \in \Om_k,
\end{eqnarray*}
 such that $\lim_{k\to\infty}\ell_k$ existed in $\g^*$. Hence the sequences
$$\{x_k\}_k,\ \{y_k\}_k,\ \{m_k\}_k,\ \{n_k\}_k,\ \{\al e^{-r_k}\}_k,\
\{\beta e^{s_k}\}_k,\
\{\al r_k-\beta s_k +\tau_k\}_k$$
would converge. Since $\tau_k$ tends to $+\infty$, it followed that
$\al r_k-\beta s_k $ would tend to $-\infty$ for  $n\to\infty.$
But $\{\al r_k\}_k$ and $\{-\beta s_k\}_k$ cannot tend to $-\infty$ for
$n\to \infty,$
 since then the sequences $\{\al e^{-r_k}\}_k$ and $\{\beta e^{s_k}\}_k$
would diverge.

This contradiction shows that $\Om_k$ tends to infinity in the orbit space
as $n$ tends
to infinity. In particular, we see that in this example
 $$\lim_{k\to \infty}||\pi_{\ell+\tau_k T^*}(f)||=0,$$
for every $f\in L^1(G)$.

\bigskip
\hfill Q.E.D.



\vskip1cm

\noindent\emph{Institute of Mathematics, Wroclaw University,
Pl. Grunnwaldzki 2/4, 50-384 Wroclaw, Poland\\ e-mail:
hebisch@math.uni.wroc.pl}
\vskip0.5cm

\noindent\emph{Universit\'e de Metz, Math\'ematiques, Ile du Saulcy,
57045 Metz Cedex, France\\ e-mail: ludwig@poncelet.univ-metz.fr}
\vskip0.5cm

\noindent\emph{Mathematisches Seminar,  C.A.-Universit\"at Kiel,
Ludewig-Meyn-Str.4, D-24098 Kiel, Germany\\ e-mail: mueller@math.uni-kiel.de}


\begin{thebibliography}{99}


\bibitem{Bernat}P.~Bernat et al.,
Repr\'esentations des groupes de Lie r\'esolubles,
{\em Dunod, Paris} 1972.

\bibitem{Boidol}J.~Boidol,
*-Regularity  of exponential Lie groups,
{\em Invent. math.} 56 (1980), 31--238.


\bibitem{ChristMueller}M.~Christ, D.~M\"uller,
 On $L^p$ spectral multipliers for a solvable Lie group,
{\em Geom. and Funct. Anal.} 6 (1996), 860--876.

\bibitem{CoifmanWeiss}R.~R.~Coifman, G.~Weiss,
Transference methods in analysis,
{\em CBMS Regional Conference Ser., A.M.S.} 31, (1977).

\bibitem{Dixmier} J.~Dixmier, Les $C^*$-alg\'ebres et leurs repr\'esentations,
{\em  Gauthier-Villars, Paris} 1964.

\bibitem{Eymard} P.~Eymard, Alg\`ebres $A_p$ et convoluteurs de $L^p$,
{\em S\'eminaire Bourbaki 1969/70, Lecture Notes in Math. } 180, {\em
Springer, Berlin 1971},
pp. 55--77.

\bibitem{Figa} A.~Fig\`a-Talamanca, Multipliers of $p$-integrable functions,
{\em Bull. Amer. Math. Soc.} 70 (1964), 666--669.

\bibitem{Herzt} C.~S.~Herz, The theory of $p$-spaces with application 
to convolution operators, {Trans. Amer. Math. Soc.} 154 (1971), 69--82.

\bibitem{Herza} C.~S.~Herz, 
Harmonic synthesis for subgroups, 
{\em Ann. Inst. Fourier} XXIII.3 (1973), 69--82.

\bibitem{Herz} C.~S.~Herz, Une g\'eneralisation de la notion de
transform\'ee de
Fourier-Stieltjes,
{\em Ann. Inst. Fourier} XXIV.3 (1974), 145--157.

\bibitem{Hoermander} L.~H\"ormander, Hypoelliptic second-order differential
equations,
{\em Acta Math.} 119 (1967), 147--171.

\bibitem{Kato} T.~Kato, Perturbation theory for linear operators,
{\em Springer, New York} 1966.

\bibitem{Krasnoselskii} M.~A.~Krasnoselskii, On a theorem of M.~Riesz,
{\em Dokl. Akad. Nauk} 131 (1959), 246--248.

\bibitem{Leptin} H.~Leptin, Sur l'alg\`ebre de Fourier d'un groupe
localement compact,
{\em C. R. Acad. Sci. Paris} 266 A (1968), 1180--1182.

\bibitem{LeptinLudwig} H.~Leptin, J.~Ludwig, Unitary representation theory of
exponential Lie groups,
{\em De Gruyter Expositions in Mathematics} 18, 1994.

\bibitem{Ludwig-Mueller} J.~Ludwig, D.~M\"uller, Sub-Laplacians of
holomorphic $L^p$-type
on rank one $AN$-groups and related solvable groups, {\em J. of Funct. Anal.}
170 (2000), 366--427.

\bibitem{MarcZygmund} J. Marcinkiewicz, A. Zygmund, Quelques in\'egalit\'es 
pour les ope\`rations lin\`eaires, {\em Fund. Math.} 32 (1939), 115--121.

\bibitem{Nelson} E.~Nelson, W.~F.~Stinespring, Representations of elliptic
operators in an enveloping algebra,
{\em  Amer. J. Math.} 81  (1959), 547--560.

\bibitem{Pier}J.-P.~Pier, Amenable Banach algebras,
{\em Pitman Research Notes in Mathematics} 172, 1988.

\bibitem{Poguntke}D.~Poguntke, Aufl\"osbare Liesche Gruppen mit symmetrischen
$L^1$-Algebren,
 {\em J. f\"ur die Reine und Angew. Math.} 358 (1985), 20--42.

\bibitem{Poguntke-oral}D.~Poguntke, oral communication.

\bibitem{Stein} E.~M.~Stein, Singular integrals and differentiability
properties
of functions, {\em Princeton Univ. Press} 1970.



\end{thebibliography}
\end{document}